\documentclass[12pt]{article}
\usepackage{amsfonts}
\input{epsf.sty}
\usepackage[dvips]{graphicx}
\usepackage{latexsym,amsmath,amsfonts,amscd, amsthm}
\usepackage{epsfig}
\usepackage{changebar}
\usepackage{pstricks}
\usepackage{pst-plot}
\usepackage{multirow}
\usepackage{subfigure}
\usepackage{placeins}
\usepackage{amssymb}

\usepackage{mathrsfs}
\usepackage{amssymb,amsthm,hyperref,amsmath,graphicx,bm,color,booktabs}
\usepackage{bm}
\usepackage{comment}

\numberwithin{equation}{section} 

\newtheorem{rem}{Remark}[section]
\theoremstyle{definition}

\newcommand{\brac}[1]{\left(#1\right)}

\newcommand{\order}{\mathcal{O}}

\newcommand{\ie}{{\it{i.e.}}}
\newcommand{\eg}{{\it{e.g.\ }}}

\newcommand{\diff}{\mbox{d}}

\allowdisplaybreaks

\begin{document}

\baselineskip=2pc

\begin{center}
{\bf \large Boundary treatment of implicit-explicit Runge-Kutta method for hyperbolic systems with source terms}
\end{center}

\centerline{
Weifeng Zhao\footnote{Department of Applied Mathematics, University of Science and Technology Beijing, Beijing 100083, China.
E-mail: wfzhao@ustb.edu.cn. Research supported by the National Natural Science Foundation of China
(No. 11801030, No. 11861131004).}
,
Juntao Huang\footnote{Department of Mathematics, Michigan State University, East Lansing, MI 48824, USA.
E-mail: huangj75@msu.edu.}
}

\vspace{.2in}

\centerline{\bf Abstract}


In this paper, we develop a high order finite difference boundary treatment method for the implicit-explicit (IMEX) Runge-Kutta (RK) schemes solving hyperbolic systems with possibly stiff source terms on a Cartesian mesh. The main challenge is how to obtain the solutions at ghost points resulting from the wide stencil of the interior high order scheme. We address this problem by combining the idea of using the RK schemes at the boundary and an inverse Lax-Wendroff procedure. The former preserves the accuracy of the RK schemes and the latter guarantees the stability. Our method is different from the widely used approach for the explicit RK schemes by imposing boundary conditions at intermediate stages \cite{Carpenter1995}, which could not be derived for the IMEX schemes. In addition, the intermediate boundary conditions are only available for explicit RK schemes up to third order while our method applies to arbitrary order IMEX and explicit RK schemes. Moreover, the present boundary treatment method may be adapted to IMEX RK schemes solving many other partial differential equations. For a specific third-order IMEX scheme, we demonstrate the good stability and third-order accuracy of our boundary treatment through both 1D examples and 2D reactive Euler equations.
\bigskip

\vfill


\newpage

\section{Introduction}
\setcounter{equation}{0}
\setcounter{figure}{0}
\setcounter{table}{0}

Implicit-explicit (IMEX) Runge-Kutta (RK) schemes \cite{ascher1997} are a powerful tool for the time discretization of many partial differential equations (PDEs), such as convection-diffusion equations \cite{Wang2015,lu2016diffusion,Wang2018}, hyperbolic systems with stiff sources \cite{pareschi2000,pareschi2005,BR2009,BPR2013,BPR2017,BQR2018}, and the Boltzmann equation \cite{FY2013,Hu2017,Hu2018,Jin2018}.
They achieve good stability by implicitly treating the stiff terms in the PDEs and high order accuracy by using several immediate stages.
Though powerful and efficient, the IMEX RK schemes for the time-dependent PDEs are not without drawbacks.
A long-standing problem is the lack of appropriate boundary treatment method \cite{ascher1997}, which is fundamental for practical applications of IMEX RK schemes.

For the explicit RK schemes solving hyperbolic equations, a popular boundary treatment method is to impose consistent boundary conditions for each intermediate stage \cite{Carpenter1995}. However, these boundary conditions are derived only for explicit RK schemes up to third order. Though it is remedied for fourth order schemes in \cite{abarbanel1996,pathria1997}, the methods therein only apply to one-dimensional (1D) scalar equations or systems with all characteristics flowing into the domain at the boundary.
Moreover, the intermediate-stage boundary conditions are not available for IMEX RK schemes since the implicit part would result in nonlinear equations which cannot be analytically solved in general.
On the other hand, boundary conditions for intermediate stages are derived and analyzed in \cite{alonso2002NumMath,alonso2002ANM,alonso2005,alonso2016,alonso2017} for implicit RK schemes solving linear scalar equations.
At this point, we remark that another limitation of all the boundary treatments mentioned above \cite{Carpenter1995,abarbanel1996,pathria1997,alonso2002NumMath,alonso2002ANM,alonso2005,alonso2016,alonso2017} is that they only deal with 1D case where the boundary is located at a grid point, restricting their applications for 2D problems on complex domains.

Perhaps due to the difficulties of extending the above boundary treatments of explicit and implicit schemes to the IMEX case,
there are rare works on that of IMEX RK schemes. To our knowledge, the only literature on this topic is \cite{Wang2018}.
In that reference, a strategy of imposing intermediate boundary conditions is presented for one-dimensional (1D) linear convection-diffusion equations. There is a similar restriction that it does not apply to general nonlinear equations.

In this paper, we propose a boundary treatment method for the IMEX RK schemes solving hyperbolic systems with possibly stiff source terms.
High order finite difference schemes on a Cartesian mesh are used for the space discretization. Due to the wide stencil of high order schemes, ghost points are needed near the boundary where the numerical stencil is partially outside of the computational domain.  Then the main task of boundary treatment is to obtain the solutions at the ghost points. We address this problem by combining the idea of using the RK schemes at the boundary and an inverse Lax-Wendroff (ILW) procedure \cite{Tan2010,TWSN2012}. The ILW procedure is employed in \cite{Tan2010,TWSN2012} as a technique of boundary treatment for hyperbolic systems  by using the PDEs to convert spatial derivatives of inflow boundaries to time derivatives, while those of the outflow boundaries are approximated with the solutions at interior points. Nevertheless, the boundary treatment in \cite{Tan2010,TWSN2012} relies on the intermediate-stage boundary conditions \cite{Carpenter1995}, which are not available for IMEX schemes as aforementioned. To overcome this difficulty, we propose to use the RK schemes at the boundary directly, instead of imposing intermediate-stage boundary conditions. Combining this idea and the ILW procedure, we can obtain the solutions as well as their derivatives at the boundary and then the solutions of the ghosts points are computed by the Taylor expansion at the boundary. In this way, our boundary treatment not only preserves the accuracy of the IMEX RK schemes but also processes good stability. In addition, it allows the boundary to be arbitrarily located.
While our method applies to general IMEX RK schemes, we demonstrate its effectiveness with a specific third-order scheme in \cite{ascher1997} by fixing the CFL number as 0.8. The third-order accuracy of the method is verified for both 1D and 2D cases, and the numerical results of the 2D reactive Euler equations are comparable to those obtained by the third-order positivity-preserving schemes with reflective boundary conditions \cite{Huang2019jsc}.

It is interesting to note that for 1D scalar equations we show that our method applied to third-order IMEX RK scheme can yield the widely used intermediate-stage boundary conditions in \cite{Carpenter1995}, indicating the correctness of our boundary treatment for IMEX RK schemes.
Since the IMEX RK schemes include explicit schemes as special cases, the present boundary treatment also applies to explicit RK schemes. It does not rely on intermediate-stage boundary conditions and thus is expected to be valid for explicit RK schemes of order higher than three. This will be investigated in our next work. Moreover, our method may be adapted to IMEX RK schemes solving the other partial differential equations, \eg  convection-diffusion equations \cite{Wang2015,lu2016diffusion,Wang2018} and the Boltzmann equation \cite{FY2013,Hu2017,Hu2018,Jin2018}.  Many unresolved issues for the boundary treatment of RK schemes can be explored based on the idea of this work.

This paper is organized as follows. In Section 2, we introduce the IMEX RK schemes for hyperbolic systems with source terms. In Section 3, we use the 1D systems to illustrate our idea of boundary treatment. In Section 4, we show that for 1D scalar equations our boundary treatment can yield the widely used intermediate-stage boundary conditions. We extend our boundary treatment to the 2D reactive Euler equations in Section 5.
Section 6 is devoted to numerical validations of our method. Finally, some conclusions and remarks are given in Section 7.

\section{IMEX RK schemes for hyperbolic systems with source terms}

Consider a one-dimensional nonlinear hyperbolic system with source term
\begin{equation}\label{21}
\partial_t U + \partial_x F(U) = Q(U),
\end{equation}
where $0 \leq x \leq 1$ and $U, F(U), Q(U) \in \mathbb R^M$.
In the domain $0 \leq x \leq 1$, appropriate boundary conditions and initial data $U(0,x)$ are required to close the above system.
We assume that the Jacobian matrix $F_{U}(U(t,0))$ always has $p$ positive eigenvalues and thus $p$ independent relations among incoming and outgoing modes are given at the left boundary $x=0$:
\begin{equation}
B(U(t,0), t) = 0.   \label{23a}
\end{equation}
Similarly, the boundary condition at $x=1$ are also imposed properly.

A general IMEX RK scheme for the hyperbolic system \eqref{21} is \cite{BR2009}
\begin{subequations}\label{24}
\begin{align}
& U^{(i)} = U^{n} - \Delta t \sum_{j=1}^{i-1} \tilde a_{ij} \partial_x F(U^{(j)}) + \Delta t \sum_{j=1}^i a_{ij}Q(U^{(j)})  ,
    \quad i=1,2,\ldots,s, \label{24a}\\
& U^{n+1} = U^{n} - \Delta t \sum_{i=1}^{s} \tilde \omega_{i} \partial_x F(U^{(i)}) + \Delta t \sum_{i=1}^s \omega_{i}Q(U^{(i)}).
   \label{24b}
\end{align}
\end{subequations}
Here the $s \times s$ matrices $\tilde A = (\tilde a_{ij})$ and $A = (a_{ij})$ satisfy $\tilde a_{ij}=0$ for $j \geq i$ and $a_{ij}=0$ for $j>i$
so that the scheme is explicit for the convection part and implicit for the source term.
Along with the coefficient vectors $\tilde w = (\tilde w_1, \tilde w_2, \ldots, \tilde w_s)^T$ and $w = (w_1, w_2, \ldots, w_s)^T$, the scheme can be represented by a double Butcher tableau
\begin{center}
\begin{tabular}{p{0.3cm} | p{0.3cm}}
$\tilde c$  &  $\tilde A$     \\
\hline
  & $\tilde w^T$
\end{tabular}
\qquad
\begin{tabular}{p{0.3cm} | p{0.3cm}}
$c$  &  $A$     \\
\hline
  & $w^T$
\end{tabular}
\end{center}
where $\tilde c = (\tilde c_1, \tilde c_2, \ldots, \tilde c_s)^T$, $c = (c_1, c_2, \ldots, c_s)^T$ are defined as
\begin{equation}
\tilde c_i = \sum_{j=1}^{i-1} \tilde a_{ij}, \quad
c_i = \sum_{j=1}^{i} a_{ij}.
\end{equation}

For the space discretization in \eqref{24}, we use the finite difference WENO scheme on a Cartesian mesh \cite{JS1996}, which needs a stencil of several points in each direction. As a consequence, ghost points are needed near the boundary $\partial \Omega$ where the numerical stencil is partially outside of $\Omega$. For explicit RK schemes, solutions at ghost points can be determined as in \cite{Tan2010,TWSN2012} with the aid of an ILW procedure and the boundary condition of each stage constructed in \cite{Carpenter1995}. However, this approach does not apply to IMEX schemes since the boundary conditions of intermediate stages could not be analytically derived for general nonlinear systems. Thus boundary treatment of IMEX RK schemes remains a challenge and we address this problem in the following.

\section{Boundary treatment for 1D system}

In this section, we use the one-dimensional system \eqref{21} to illustrate our idea of boundary treatment for the IMEX RK scheme \eqref{24}, and then extend it to the two-dimensional case in Section \ref{sec5}. We employ the third-order IMEX RK scheme for concreteness, while our method applies to general IMEX schemes. Additionally, we use the third-order finite difference WENO scheme in space \cite{Liu1994weno3} and take $\Delta t=\order(\Delta x)$.

\subsection{Computation of solutions at ghost points}\label{sec31}

We focus on the left boundary $x=0$ of the problem \eqref{21} and the method can be similarly applied to the right boundary.
Following the notations in \cite{TWSN2012}, we discretize the interval $(0,1)$ by a uniform mesh
\begin{equation}\label{3-1}
\frac{\Delta x}{2} = x_0 < x_1 < \cdots <x_N = 1- \frac{\Delta x}{2}
\end{equation}
and set $x_j, j=-1,-2$ as two ghost points near the left boundary $x=0$ (note that the boundary indeed can be arbitrarily located).
Denote by $U_j^n$ the numerical solution of $U$ at position $x_j$ and time $t_n$.
Assume that the interior solutions $U_j, j=0,1,2,\ldots,N$ have been updated from time level $t_{n-1}$ to time level $t_n$.
As in \cite{TWSN2012}, we use a third-order Taylor approximation to construct the values at the ghost points $x_j$ for $j=-1,-2$, namely,
\begin{equation}\label{25}
U_j^n = \sum_{k=0}^{2} \frac{x_j^k}{k!} U^{n,(k)} , \quad j=-1,-2.
\end{equation}
Here $U^{n,(k)}$ denotes a $(3-k)$-th order approximation of the spatial derivative at the boundary point $\frac{\partial^k U}{\partial x^k}\big|_{x=0, t=t_n}$.
With this formula, $U_j^n$ at the ghost points can be obtained once $U^{n,(k)}, k=0,1,2$ are provided, which is well treated in \cite{Tan2010,TWSN2012} (see also the next subsection for details).
Then we can obtain  $U_j^{(1)}$ at the interior points by the RK scheme \eqref{24a}. The next task is to compute $U_j^{(1)}, j=-1,-2$.
We propose an approach to do this as follows.

Similar to the above procedure, we compute $U_j^{(1)}$ with the third-order Taylor expansion at the boundary point $x_b=0$:
\begin{equation}\label{27}
U_j^{(1)} = \sum_{k=0}^{2} \frac{x_j^k}{k!} U^{(1),(k)},  \quad j=-1,-2,
\end{equation}
where $U^{(1),(k)}$ denotes a $(3-k)$-th order approximation of the spatial derivative $\frac{\partial^k U^{(1)}}{\partial x^k}\big|_{x=0}$.
Then we turn to compute $U^{(1),(k)}$ for $k=0,1,2$.
To this end, we apply the first stage of the RK solver \eqref{24a} for $U^{(1)}$ at the boundary point $x_b=0$:
\begin{equation}\label{28}
U^{(1)}(x_b) = U^{n}(x_b) - \tilde a_{10}  \Delta t \partial_x F(U^n(x_b)) + \Delta t a_{10} Q(U^{(1)}(x_b))  .
\end{equation}
Notice that $U^{n}(x_b) = U^{n,(0)} + \order(\Delta x^3)$ and
$\partial_x F(U^n(x_b)) = F_U(U^n(x_b)) \partial_x U^n(x_b) = F_U(U^{n,(0)}) U^{n,(1)} + \order(\Delta x^2)$ are already known.
Substituting these into \eqref{28} and solving the algebra equation, we can obtain an approximation of  $U^{(1)}(x_b)$ and denote it by $U^{(1),(0)}$. In addition, the error between $U^{(1),(0)}$ and  $U^{(1)}(x_b)$ defined by \eqref{28} is $\order(\Delta x^3)$.


Furthermore, taking derivatives with respect to $x$ on both sides of \eqref{28} yields
\begin{equation}\label{29}
\frac{\partial U^{(1)}}{\partial x}\big|_{x= x_b}
= \partial_x U^{n} (x_b)
 -\tilde a_{10} \Delta t \partial_{xx} F(U^n(x_b)) + \Delta t a_{10} Q_U(U^{(1)}(x_b)) \frac{\partial U^{(1)}}{\partial x}\big|_{x=x_b} .
\end{equation}
Here $\partial_x U^{n} (x_b) = U^{n,(1)} + \order(\Delta x^2)$ and
$\partial_{xx} F(U^n(x_b)) = F_{UU}(U^n(x_b)) \partial_x U^n(x_b) \partial_x U^n(x_b) + F_U(U^n(x_b)) \partial_{xx} U^n(x_b)
= F_{UU}(U^{n,(0)}) U^{n,(1)} U^{n,(1)} + F_U(U^{n,(0)}) U^{n,(2)} +  \order(\Delta x)$.
With these approximations, $\frac{\partial U^{(1)}}{\partial x}\big|_{x= x_b} $ can be solved out from \eqref{29} and the resulting solution  $U^{(1),(1)}$ is a second-order approximation of $\frac{\partial U^{(1)}}{\partial x}\big|_{x= x_b} $.

By taking higher order derivatives on both sides of \eqref{28}, one can also compute $\frac{\partial ^k U^{(1)}}{\partial x^k}\big|_{x= x_b}$ for $k=2$. However, this procedure is quite complicated as it involves Jacobian of Jacobian. Here we simply approximate $\frac{\partial ^k U^{(1)}}{\partial x^k}\big|_{x= x_b} $ for $k=2$ by using the $(3-k)$-th order WENO type extrapolation with $U^{(1)}$ at interior points in \cite{TWSN2012}, which will be introduced in subsection \ref{sec33}. In this way, we obtain $U^{(1),(k)}$ for $k=0,1,2$ with accuracy of order $(3-k)$. Then $U_j^{(1)}$ for $j=-1,-2$ can be computed by \eqref{27} with third-order accuracy. Having $U^{(1)}$ at the ghost points, we can then evolve from $U^{(1)}$ to $U^{(2)}$ using the interior difference scheme.

Repeating the same procedure for each $U^{(i)}, 2\le i\le s$, we can compute the solution at the ghost points in the $i$-th intermediate stage. Then we can update the solution at all interior points in the $(i+1)$-th stage. Finally, we obtain $U^{n+1}$, \ie, the solution at the end of a complete RK cycle.

\subsection{Computation of $U^{n,(k)}$ at the boundary}\label{sec32}

For the sake of completeness, we provide the method in \cite{TWSN2012} for computing $U^{n,(k)}, k=0,1,2$, \ie, the $(3-k)$-th order approximation of $\frac{\partial^k U^n}{\partial x^k}$ at the boundary $x=0$.

\subsubsection{$k=0$}


We first do a local characteristic decomposition to determine the inflow and outflow boundary conditions as in \cite{TWSN2012}. Denote the Jacobian matrix of the flux evaluated at $x=x_0$ by
$$A(U_0^n) = \partial_U F(U)\big|_{U=U_0^n}$$
and assume that it has $p$ positive eigenvalues $\lambda_1, \lambda_2, \ldots, \lambda_p$ and $(M-p)$ negative eigenvalues $\lambda_{p+1}, \lambda_{p+2}, \ldots, \lambda_M$ with $l_1, l_2, ...l_p$ and $l_{p+1}, l_{p+2}, ...l_M$ the corresponding left eigenvectors, respectively.
Define by $V_{j,m}$ the $m$-th component of local characteristic variable at grid point $x_j,j=0,1,2$, \ie,
\begin{equation}\label{310}
V_{j,m} = l_m U_j^n, \quad  m=1,2,\ldots,M, \quad j=0,1,2.
\end{equation}
We extrapolate the outgoing characteristic variable $V_{j,m}, m=p+1,p+2,\ldots,M$ to the boundary $x_b$ with the WENO type extrapolation \cite{TWSN2012} (see also subsection \ref{sec33}), and denote the extrapolated $k$-th order derivative by
\begin{equation}\label{311}
V_{x_b, m}^{*(k)}, \quad k=0,1,2.
\end{equation}
With $V_{x_b, m}^{*(0)}$ and the boundary condition \eqref{23a}, $U^{n,(0)}$ at the boundary can be solved out from the following equation
\begin{equation}\label{312_0}
\begin{split}
& l_m U^{n,(0)} = V_{x_b, m}^{*(0)},  \quad m = p+1,p+2,\ldots,M, \\
& B( U^{n,(0)}, t_n ) = 0.
\end{split}
\end{equation}

\subsubsection{$k=1$}\label{sec322}

Having $U^{n,(0)}$, we proceed to compute $U^{n,(1)}$ with the ILW procedure proposed in \cite{Tan2010,TWSN2012}.
To do this, we take derivative with respect to $t$ for the boundary condition \eqref{24}
\begin{equation*}
B_U( U^{n,(0)}, t_n  ) \partial_t U^{n,(0)} + B_t( U^{n,(0)}, t_n  ) = 0,
\end{equation*}
which can be written as
\begin{equation*}
B_U( U^{n,(0)}, t_n  ) \partial_t U^{n,(0)} = g(U^{n,(0)} , t_n )
\end{equation*}
with $g(U^{n,(0)} , t_n ) := -B_t( U^{n,(0)}, t_n  )$. On the other hand, multiplying the equation \eqref{21} with $B_U( U^{n,(0)}, t_n  )$ from the left yields
$$
B_U( U^{n,(0)}, t_n  ) \partial_t U^{n,(0)}  +  B_U( U^{n,(0)}, t_n  ) A( U^{n,(0)} ) U^{n,(1)} = B_U( U^{n,(0)}, t_n  )Q(U^{n,(0)}).
$$
With the above two equations and $V_{x_b, m}^{*(1)}$ obtained by the WENO type extrapolation, $U^{n,(1)}$ can be solved out from
\begin{subequations}\label{312}
\begin{align}
& l_m U^{n,(1)} = V_{x_b, m}^{*(1)},  \quad m = p+1,p+2,\ldots,M,  \label{312a}\\
& B_U( U^{n,(0)}, t_n  ) A( U^{n,(0)} ) U^{n,(1)} = B_U( U^{n,(0)}, t_n  )Q(U^{n,(0)}) - g(U^{n,(0)} , t_n ). \label{312b}
\end{align}
\end{subequations}
This is the ILW procedure in \cite{Tan2010,TWSN2012}.


\subsubsection{$k=2$}
Following \cite{TWSN2012}, we approximate $\frac{\partial^2 U^n}{\partial x^2}|_{x=0}$ with the WENO type extrapolation. Specifically, since we have obtained the characteristic variables at grid points near the boundary in \eqref{310}, the WENO type extrapolation is employed based on these characteristic variables to compute the second-order derivative $V_{x_b, m}^{*(2)}$ for each $m=1,2,\cdots, M$. Then the approximation of $\frac{\partial^2 U^n}{\partial x^2}|_{x=0}$ is given by
\begin{equation}\label{3}
U^{n,(2)} = L^{-1} V_{x_b, m}^{*(2)},
\end{equation}
where $L$ is local left eigenvector matrix composed of $l_1, l_2, \cdots, l_M$.  Similarly, $\frac{\partial^k U^n}{\partial x^k}|_{x=0}, k>2$ are also computed with the WENO type extrapolation if necessary for higher order boundary treatments.


\subsection{WENO type extrapolation}\label{sec33}

Finally, we show how to obtain a $(3-k)$-th order approximation of $\frac{\partial^k V}{\partial x^k}|_{x=0}$, denoted by $V^{*(k)}, k=0,1,2$, with the stencil $S_3=\{ x_0, x_1, x_2\}$, which is introduced in \cite{Tan2010}. When $V$ is smooth on $S_3$, we compute $V^{*(k)}$ as
\begin{equation*}
V^{*(k)} = \sum_{r=0}^2 d_r \frac{\diff^k p_r(x)}{\diff x^k} \big|_{x=0},
\end{equation*}
where $d_0 = \Delta x^2, d_1 = \Delta x, d_2 = 1- \Delta x- \Delta x^2$ and $p_r(x)$ is the $r$-th order Lagrange polynomial constructed with $S_3$.
When there is a discontinuity in the stencil $S_3$, the WENO type extrapolation is applied \cite{Tan2010}:
\begin{equation}\label{weno}
V^{*(k)} = \sum_{r=0}^2 \omega_r \frac{\diff^k p_r(x)}{\diff x^k} \big|_{x=0},
\end{equation}
where $\omega_r$ are nonlinear weights depending on the value of $V_j$. Similar to the classical WENO scheme, the nonlinear weights are given by \cite{Tan2010}
\begin{equation*}
\omega_r = \frac{\alpha_r}{\sum_{s=0}^2 \alpha_s}, \quad \alpha_r = \frac{d_r}{(\varepsilon + \beta_r)^2}.
\end{equation*}
Here $\varepsilon = 10^{-6}$ and $\beta_r$ are the smoothness indicators determined by
\begin{equation*}
\beta_0 = \Delta x^2, \quad
\beta_r = \sum_{l=1}^r \int_{x_{-1}}^{x_0} \Delta x^{2l-1}
\left(   \frac{\diff^l}{\diff x^l} p_r(x)  \right)^2 \diff x, \quad r=1,2.
\end{equation*}
The explicit expressions of $\beta_1$ and $\beta_2$ are given in \cite{Tan2010}. When $V$ is smooth on $S_3$, it is shown in \cite{Tan2010} that
\begin{equation*}
\omega_0 = \order(\Delta x^2), \quad
\omega_1 = \order(\Delta x), \quad
\omega_2 = 1-\omega_0-\omega_1
\end{equation*}
so that the WENO type extrapolation \eqref{weno} is $(3-k)$-th order accurate.

In the end of this section, we summarize our third-order boundary treatment for the hyperbolic system \eqref{21} as follows:

{\bf Step 1:} Compute $U^{n,(k)}$, the $(3-k)$-th order approximation of $\frac{\partial^k U^n}{\partial x^k}$ at the boundary $x=0$ for $k=0,1,2$ as in subsection \ref{sec32}. Then impose $U^n_j$ at the ghost points $j=-1,-2$ with the Taylor expansion \eqref{25}. This step is the same as that in \cite{TWSN2012} for explicit RK schemes. With $U^n$ at the ghost points, we can evolve from $U^n$ to $U^{(1)}$ with the interior difference scheme.

{\bf Step 2:} Compute $U^{(1),(k)}$, the $(3-k)$-th order approximation of $\frac{\partial^k U^{(1)}}{\partial x^k}$ at the boundary $x=0$ for $k=0,1,2$ as in subsection \ref{sec31}. Specifically, $U^{(1),(0)}$ and $U^{(1),(1)}$ are solved out from the RK solver \eqref{28} and \eqref{29}, respectively, and $U^{(1),(2)}$ is approximated with the WENO type extrapolation in subsection \ref{sec33}. Then impose $U^{(1)}_j$ at the ghost points $j=-1,-2$ with the Taylor expansion \eqref{27}. Having $U^{(1)}$ at the ghost points, we can update $U^{(2)}$ at all interior points.

{\bf Step 3:} Repeat the same procedure as in Step 2 for each $U^{(i)}, 2\le i\le s$ to compute the solution at the ghost points in the $i$-th intermediate stage. Then we can update the solution at all interior points in the $(i+1)$-th stage. Finally, we obtain $U^{n+1}$, \ie, the solution at the end of a complete RK cycle.

\begin{rem}
Combining the idea of using the RK schemes of intermediate solutions at the boundary and the ILW procedure is the key of our boundary treatment. The former preserves the accuracy of the RK method and the latter guarantees the stability.
This is different from the widely used approach for the explicit RK scheme by imposing boundary conditions for intermediate solutions \cite{Carpenter1995}, which are unavailable for the IMEX schemes. In addition, the intermediate boundary conditions are only available for explicit RK schemes up to third order while our method applies to general IMEX RK scheme.
\end{rem}

\begin{rem}
Though the idea of using the RK schemes of intermediate solutions has been stressed in \cite{pathria1997}, the method therein is only valid for
1D scalar equation or 1D system with all characteristics flowing into the domain at the boundary. Another limitation of the method in \cite{pathria1997} is that the boundary has to be located at a grid point, which makes it difficult to be extended to 2D case.
\end{rem}

\begin{rem}
Since the IMEX RK schemes include explicit schemes as special cases, the present boundary treatment also applies to explicit RK schemes.
It does not need to impose intermediate boundary conditions and thus is expected to be valid for explicit RK schemes of order higher than three. This will be investigated in our next work.
\end{rem}



\section{Analysis of the boundary treatment}

In this section, we show that our boundary treatment applied to the third-order explicit RK scheme solving 1D scalar equations can yield the widely used boundary conditions for intermediate solutions in \cite{Carpenter1995}.

\subsection{Linear equation}

We consider a linear scalar equation
\begin{subequations}\label{4-0}
\begin{align}
& \partial_t u + \partial_x u = 0, \quad 0 \leq x \leq 1, \quad t >0, \\
& u(t,0) = g(t)
\end{align}
\end{subequations}
and solve it with the third-order explicit strong-stability-preserving RK scheme \cite{Gottlieb2001sspReview}
\begin{subequations}\label{4-00}
\begin{align}
& u^{(1)} = u^n + \Delta t \mathcal{L}( u^n  ), \label{4-00a} \\
& u^{(2)} = \frac{3}{4}u^n + \frac{1}{4}u^{(1)}
               + \frac{1}{4} \Delta t \mathcal{L}( u^{(1)}  ), \label{4-00b}\\
& u^{n+1} = \frac{1}{3}u^n + \frac{2}{3}u^{(2)}
               + \frac{2}{3} \Delta t \mathcal{L}( u^{(2)}  ), \label{4-00c}
\end{align}
\end{subequations}
where $\mathcal{L}( u ) = - \partial_x u$. We focus on the boundary treatment at the left boundary $x_b=0$.

According to the boundary treatment in the previous section, we first need to compute $u^{n}$ and its first- and second-order derivatives at the boundary point $x=x_b$, denoted by $u^n_b$ and $u^{n,(k)}_b, k=1,2$, respectively.
For the above scalar equation, $u^{n}_b$ is simply given by the boundary condition:
\begin{equation}\label{4-1}
u^{n}_b = g(t_n).
\end{equation}
Then according to subsection \ref{sec322}, the ILW procedure is employed to compute the first-order derivative $u^{n,(1)}_b$:
\begin{equation}\label{4-2}
u^{n,(1)}_b = \partial_x u( t_n, x_b ) = - \partial_t u( t_n, x_b ) = - g^{\prime}(t_n).
\end{equation}
Similarly, the second-order derivative $u^{n,(2)}_b$ computed by the ILW procedure is 
\begin{equation}\label{4-3}
u^{n,(2)}_b = g^{\prime \prime}(t_n),
\end{equation}
and by the third-order WENO extrapolation is 
\begin{equation}\label{4-4}
u^{n,(2)}_b = g^{\prime \prime}(t_n) + \order(\Delta x).
\end{equation}
Here the solution is assumed to be smooth near the boundary for \eqref{4-4}.

Next we compute $u^{(1)}_b$ and its derivatives $u^{(1),(k)}_b, k=1,2$. In our method, the RK scheme for $u^{(1)}$ is applied to determine its value at the boundary, \ie,
\begin{equation}\label{4-5}
u^{(1)}_b = u^n_b - \partial_x u^n_b = g(t_n) + \Delta t g^{\prime}(t_n),
\end{equation}
where \eqref{4-1} and \eqref{4-2} have been used. As shown in \eqref{29}, we take derivatives with respect to $x$ on both sides of
\eqref{4-00a} and evaluate at $x_0$ to obtain
\begin{equation*}
\partial_x u^{(1)}_b = \partial_x u^n_b - \Delta t \partial_{xx} u^n_b,
\end{equation*}
and thereby
\begin{equation}\label{4-6}
u^{(1),(1)}_b = - g^{\prime}(t_n) - \Delta t u^{n,(2)}_b,
\end{equation}
where \eqref{4-2} has been used. Applying the scheme \eqref{4-00b} for $u^{(2)}$ at the boundary $x_0$ and using \eqref{4-1}, \eqref{4-5} and \eqref{4-6}, we arrive at
\begin{equation}\label{4-7}
\begin{split}
u^{(2)}_b
& = \frac{3}{4}u^n_b + \frac{1}{4}u^{(1)}_b
           - \frac{1}{4} \Delta t u^{(1),(1)}_b, \\
& =  g(t_n) + \frac{1}{2} \Delta t g^{\prime}(t_n)   +  \frac{1}{4} \Delta t^2 \Delta t u^{n,(2)}_b.
\end{split}
\end{equation}
If we use \eqref{4-3} for $u^{n,(2)}_b$ obtained by the ILW procedure, then
\begin{equation}\label{4-8}
u^{(2)}_b  = g(t_n) + \frac{1}{2} \Delta t g^{\prime}(t_n)   +  \frac{1}{4} \Delta t^2 g^{\prime \prime}(t_n).
\end{equation}
In this case, the relations \eqref{4-1}, \eqref{4-5} and \eqref{4-8} are exactly the boundary conditions imposed at each stage in \cite{Carpenter1995,Tan2010}.

If the WENO extrapolation is adopted for $u^{n,(2)}_b$, then
\begin{equation}\label{4-9}
u^{(2)}_b  = g(t_n) + \frac{1}{2} \Delta t g^{\prime}(t_n)   +  \frac{1}{4} \Delta t^2 g^{\prime \prime}(t_n) + \order(\Delta t^3),
\end{equation}
which only differs from that in \eqref{4-8} with a third-order term $\order(\Delta t^3)$ and thus maintains the third-order accuracy of the boundary treatment.

\subsection{Nonlinear equation}

We further show that the above analysis also applies to the nonlinear case. To this end, we use the third-order RK scheme \eqref{4-00} for a nonlinear scalar equation
\begin{equation*}
\begin{split}
& \partial_t u + \partial_x f(u) = 0, \quad 0 \leq x \leq 1, \quad t >0, \\
& u(t,0) = g(t)
\end{split}
\end{equation*}
with $ f^{\prime}(u(0,t)) >0$, and apply our method at the left boundary $x_b=0$.

First, it is trivial that the boundary condition for $u^n$ is $u^n_b = g(t_n)$. Using the ILW procedure, we can obtain \cite{Tan2010}
\begin{equation}\label{4-10}
\begin{split}
& u^{n,(1)}_b = \partial_x u(t_n, x_b) = -\frac{g^{\prime}(t_n)}{f^{\prime}(g(t_n))} \\
& u^{n,(2)}_b = \partial_{xx} u(t_n, x_b) = \frac{f^{\prime}(g(t_n)) g^{\prime\prime}(t_n) - 2f^{\prime\prime}(g(t_n)) g^{\prime}(t_n)^2}{f^{\prime}(g(t_n))^3}.
\end{split}
\end{equation}
According to our boundary treatment, the RK scheme \eqref{4-00a} for $u^{(1)}$ is used at the boundary to determine $u^{(1)}_b$:
\begin{equation}\label{4-11}
u^{(1)}_b = u^n_b - \Delta t f^{\prime}( u^n_b  ) u_x(t_n, x_b)
= g(t_n) + \Delta t g^{\prime}(t_n).
\end{equation}
Moreover, taking derivatives with respect to $x$ on both sides of the scheme \eqref{4-00a} and evaluating at the boundary point $x=x_b$, we have
\begin{equation}\label{4-12}
\begin{split}
u^{(1),(1)}_b
&= u^{n,(1)}_b - \Delta t \brac{f^{\prime \prime}( u^n_b  ) \partial_xu^2(t_n, x_b) + f^{\prime}( u^n_b  ) \partial_{xx} u(t_n, x_b)} \\
&= -\frac{g^{\prime}(t_n)}{f^{\prime}(g(t_n))} - \Delta t f^{\prime \prime}( g(t_n) ) \left(  \frac{g^{\prime}(t_n)}{f^{\prime}(g(t_n))} \right)^2 \\
&\hspace{5mm} - \Delta t f^{\prime}( g(t_n)  ) \frac{f^{\prime}(g(t_n)) g^{\prime\prime}(t_n)
     - 2f^{\prime\prime}(g(t_n)) g^{\prime}(t_n)^2}{f^{\prime}(g(t_n))^3} \\
& = -\frac{g^{\prime}(t_n)}{f^{\prime}(g(t_n))} - \Delta t \frac{ f^{\prime}(g(t_n)) g^{\prime\prime}(t_n) - f^{\prime\prime}(g(t_n)) g^{\prime}(t_n)^2 }{f^{\prime}(g(t_n))^2}.
\end{split}
\end{equation}
Thus $u^{(2)}_b$ determined by the scheme \eqref{4-00b} is
\begin{equation}\label{4-13}
\begin{split}
u^{(2)}_b
&= \frac{3}{4} u^{n}_b + \frac{1}{4} u^{(1)}_b
   - \frac{1}{4}\Delta t f^{\prime }( u^{(1)}_b  ) u^{(1),(1)}_b \\
&=  g(t_n) + \frac{1}{4} \Delta t g^{\prime}(t_n) -  \frac{1}{4} \Delta t f^{\prime }( g(t_n) + \Delta t g^{\prime}(t_n) )\\
& \hspace{5mm} \times \brac{-\frac{g^{\prime}(t_n)}{f^{\prime}(g(t_n))}
    - \Delta t \frac{ f^{\prime}(g(t_n)) g^{\prime\prime}(t_n) - f^{\prime\prime}(g(t_n)) g^{\prime}(t_n)^2 }{f^{\prime}(g(t_n))^2} }\ \\
&=  g(t_n) + \frac{1}{4} \Delta t \left( 1+ \frac{ f^{\prime }( g(t_n) + \Delta t g^{\prime}(t_n) ) }{ f^{\prime}(g(t_n)) } \right) g^{\prime}(t_n)\\
&  \hspace{5mm} +  \frac{1}{4} \Delta t^2
  \frac{ f^{\prime }( g(t_n) + \Delta t g^{\prime}(t_n) )
        \brac{f^{\prime}(g(t_n)) g^{\prime\prime}(t_n) - f^{\prime\prime}(g(t_n)) g^{\prime}(t_n)^2} }
       {f^{\prime}(g(t_n))^2}.
\end{split}
\end{equation}
Equations \eqref{4-11} and \eqref{4-13} are the boundary conditions for $u^{(1)}$ and $u^{(2)}$, respectively.
In particular, \eqref{4-11} is the same with \eqref{4-5}, and \eqref{4-13} degenerates to \eqref{4-8} when $f(u)$ is linear.
Moreover, for nonlinear $f(u)$ we conduct the Taylor expansion
\begin{equation*}
f^{\prime }( g(t_n) + \Delta t g^{\prime}(t_n) )
= f^{\prime }( g(t_n) ) + \Delta t g^{\prime}(t_n)f^{\prime \prime }( g(t_n) ) + \order( \Delta t^2 )
\end{equation*}
and then \eqref{4-13} reduces to
\begin{equation*}
\begin{split}
u^{(2)}_b
&=  g(t_n) + \frac{1}{4} \Delta t \left( 2 +  \Delta t\frac{ g^{\prime}(t_n)f^{\prime \prime }( g(t_n) ) }{ f^{\prime}(g(t_n)) } \right) g^{\prime}(t_n)\\
&  \hspace{5mm} +  \frac{1}{4} \Delta t^2
  \frac{ f^{\prime }( g(t_n))
         \left[f^{\prime}(g(t_n)) g^{\prime\prime}(t_n) - f^{\prime\prime}(g(t_n)) g^{\prime}(t_n)^2 \right]}
       {f^{\prime}(g(t_n))^2} + \order(\Delta t^3)\\
&= g(t_n) + \frac{1}{2} \Delta t g^{\prime}(t_n) + \frac{1}{4}\Delta t^2 g^{\prime\prime}(t_n) + \order(\Delta t^3).
\end{split}
\end{equation*}
The above expression of $u^{(2)}_b$ only differs from \eqref{4-8} with $\order(\Delta t^3)$ and thus the overall accuracy is still third order.

This analysis shows the correctness of our boundary treatment when it is applied to the explicit third-order RK scheme.

\section{Boundary treatment for 2D reactive Euler equations}\label{sec5}

In this section, we extend the above idea of boundary treatment to the 2D case. For convenience of representation, we concentrate on the 2D reactive Euler equations
\begin{equation}\label{eq:euler}
\partial_t U + \partial_x F(U) + \partial_y G(U) = S(U), \quad  (x,y) \in \Omega,
\end{equation}
where
\begin{equation*}
\begin{split}
& U = (\rho, \rho u, \rho v, E, \rho Y)^T, \\
& F(U) = ( \rho u, \rho u^2 + p, \rho u v, (E+p)u, \rho u Y )^T, \\
& G(U) = (\rho v, \rho u v, \rho v^2 + p, (E+p)v, \rho v Y )^T, \\
& S(U) = (0,0,0,0,\omega)
\end{split}
\end{equation*}
and
$E=\frac{1}{2}\rho(u^2 + v^2) + \frac{p}{\gamma-1} + \rho q Y$.
Here $\rho$ is the density, $u$ and $v$ are the velocities in $x$ and $y$ directions, $E$ is the total energy, $p$ is the pressure, $Y$ is the reactant mass fraction, $q>0$ is the heat release of reaction and $\gamma$ is the specific heat ratio. The source term is assumed to be in an Arrhenius form
\begin{equation*}
\omega = - \tilde K \rho Y e^{-\tilde T / T},
\end{equation*}
where $T=p/\rho$ is the temperature, $\tilde T >0$ is the activation constant temperature and $\tilde K >0$ is a constant rate coefficient.

The IMEX RK scheme for the 2D reactive Euler equations \eqref{eq:euler} is
\begin{subequations}\label{eq:IMEX-2D}
\begin{align}
& U^{(i)} = U^{n} - \Delta t \sum_{j=1}^{i-1} \tilde a_{ij} \partial_x F(U^{(j)})
                           - \Delta t \sum_{j=1}^{i-1} \tilde a_{ij} \partial_y G(U^{(j)}) + \Delta t \sum_{j=1}^i a_{ij}S(U^{(j)})  ,
    \quad i=1,2,\ldots,s, \label{24a-2d}\\
& U^{n+1} = U^{n} - \Delta t \sum_{i=1}^{s} \tilde \omega_{i} \partial_x F(U^{(i)})
                              - \Delta t \sum_{i=1}^{s} \tilde \omega_{i} \partial_y G(U^{(i)}) + \Delta t  \sum_{i=1}^s \omega_{i}S(U^{(i)}),
   \label{24b-2d}
\end{align}
\end{subequations}
where the coefficients $\tilde a_{ij}$, $a_{ij}$, $\tilde \omega_{i}$ and $\omega_{i}$ are the same with those in \eqref{24}. Assume the solutions of all the grid points inside $\Omega$ have been updated from time level $t_{n-1}$ to $t_n$. For a ghost point $P$, we find a point $P_0=(x_0,y_0)=\bm x_0$ on the boundary $\partial \Omega$ such that the normal $\bm n (\bm x_0)$ at $P_0$ goes through $P$. Here $\bm n (\bm x_0)$ points to the exterior of $\Omega$. We set up a local coordinate system at $P_0$ by
\begin{equation}
\left(
\begin{array}{c}
\hat x  \\
\hat y
\end{array}
\right)
=
\left(
\begin{array}{cc}
\cos \theta & \sin \theta \\
-\sin \theta & \cos \theta
\end{array}
\right)
\left(
\begin{array}{c}
x  \\
y
\end{array}
\right)
= \mathsf{T}\left(
\begin{array}{c}
x  \\
y
\end{array}
\right),
\end{equation}
where $\theta$ is the angle between the normal $\bm n (\bm x_0)$ and the $x$-axis, and $\mathsf{T}$ is a rotational matrix. The $\hat x$-axis then points in the same direction as $\bm n (\bm x_0)$ and the $\hat y$-axis points in the tangent direction. In this local coordinate system, the reactive Euler equations \eqref{eq:euler} become
\begin{equation}\label{eq:euler2}
\hat U_t + F(\hat U)_{\hat x} + G(\hat U)_{\hat y} = S( \hat U)
\end{equation}
with
\begin{equation*}
 \hat U = (\rho, \rho \hat u, \rho \hat v, E, \rho Y)^T,  \quad
 \left(
\begin{array}{c}
\hat u  \\
\hat v
\end{array}
\right)
= \mathsf{T}\left(
\begin{array}{c}
u \\
v
\end{array}
\right).
\end{equation*}
The IMEX scheme in this coordinate is
\begin{subequations}\label{eq:IMEX-2D-2-hat}
\begin{align}
& \hat U^{(i)} = \hat U^{n} - \Delta t \sum_{j=1}^{i-1} \tilde a_{ij} \partial_{\hat x} F(\hat U^{(j)})
                           - \Delta t \sum_{j=1}^{i-1} \tilde a_{ij} \partial_{\hat y} G(\hat U^{(j)}) + \Delta t \sum_{j=1}^i a_{ij}S(\hat U^{(j)})  ,
    \quad i=1,2,\ldots,s, \label{24a-hat}\\
& \hat U^{n+1} = \hat U^{n} - \Delta t \sum_{i=1}^{s} \tilde \omega_{i} \partial_{\hat x} F(\hat U^{(i)})
                              - \Delta t \sum_{i=1}^{s} \tilde \omega_{i} \partial_{\hat y} G(\hat U^{(i)}) + \Delta t  \sum_{i=1}^s \omega_{i}S(\hat U^{(i)}).
   \label{24b-hat}
\end{align}
\end{subequations}

For a third-order boundary treatment, $\hat U^n$ at the ghost point $P$ is imposed by the Taylor expansion
\begin{equation}
\hat U^n( P) = \sum_{k=0}^2 \frac{\Delta^k}{k!} \hat U^{n,(k)},
\end{equation}
where $\Delta$ is the $\hat x$-coordinate of $P$ and $\hat U^{n,(k)}$ is a $(3-k)$-th order approximation of the normal derivative $\frac{\partial^k \hat U}{\partial \hat x^k} \big|_{(x,y)=\bm x_0, t=t_n}$.  With $\hat U^n$ at ghost points, $\hat U^{(1)}$ at interior points can be obtained with the scheme \eqref{24a-hat}. Since the computation of $\hat U^{n,(k)}$ is the same with that in \cite{TWSN2012} and its basic idea has been illustrated for the 1D case, we omit the details here.

The next step is to compute $\hat U^{(1)}$ at the ghost points. Similar to the 1D case, we use the Taylor expansion
\begin{equation}\label{eq:U2}
\hat U^{(1)}( P) = \sum_{k=0}^2 \frac{\Delta^k}{k!} \hat U^{(1),(k)}
\end{equation}
with  $\hat U^{(1),(k)}$ a $(3-k)$-th order approximation of $\frac{\partial^k \hat U^{(1)}}{\partial \hat x^k} \big|_{(x,y)=\bm x_0}$, and turn to determine $\hat U^{(1),(k)}$, $k=0,1,2$.
First, $\hat U^{(1),(0)}=\hat U^{(1)}(P_0)$ can be computed with scheme \eqref{24a-hat}, \ie,
\begin{equation}
\hat U^{(1)}(P_0) = \hat U^{n}(P_0) - \Delta t  \tilde a_{11} \partial_{\hat x} F(\hat U^{n})(P_0)
                           - \Delta t  \tilde a_{11} \partial_{\hat y} G(\hat U^{n})(P_0) + \Delta t  a_{11}S(\hat U^{(1)}) (P_0).
\end{equation}
Furthermore, taking derivatives with respect to $\hat x$ and $\hat y$ on both sides of the above equation, we obtain
\begin{equation}
\begin{split}
& \partial_{\hat x} \hat U^{(1)}(P_0) =  \partial_{\hat x}\hat U^{n}(P_0)
                           - \Delta t  \tilde a_{11} \partial_{\hat x\hat x} F(\hat U^{n})(P_0)
                           - \Delta t  \tilde a_{11} \partial_{\hat x\hat y} G(\hat U^{n})(P_0)
                           + \Delta t a_{11}S_{U}(\hat U^{(1)})  \partial_{\hat x}\hat U^{(1)}(P_0) , \\
&  \partial_{\hat y} \hat U^{(1)}(P_0) =  \partial_{\hat y}\hat U^{n}(P_0)
                           - \Delta t  \tilde a_{11} \partial_{\hat x\hat y} F(\hat U^{n})(P_0)
                           - \Delta t  \tilde a_{11} \partial_{\hat y\hat y} G(\hat U^{n})(P_0)
                           + \Delta t a_{11}S_{U}(\hat U^{(1)})  \partial_{\hat y}\hat U^{(1)}(P_0).
\end{split}
\end{equation}
Since $\hat U^{(1)}(P_0)$ has already been obtained, and $ \partial_{\hat x\hat x} F(\hat U^{n})(P_0) , \partial_{\hat x\hat y} G(\hat U^{n})(P_0), \partial_{\hat x\hat y} F(\hat U^{n})(P_0) $ and $\partial_{\hat y\hat y} G(\hat U^{n})(P_0) $ can be approximated with $\hat U^{n}$ at interior grid points by the 2D WENO type extrapolation \cite{TWSN2012}, then $ \partial_{\hat x} \hat U^{(1)}(P_0)$ and $ \partial_{\hat y} \hat U^{(1)}(P_0)$ can be solved out from the above equation. High order derivatives, $\partial_{\hat x \hat x} \hat U^{(1)}(P_0) $, $\partial_{\hat x \hat y} \hat U^{(1)}(P_0) $ and $\partial_{\hat y \hat y} \hat U^{(1)}(P_0) $, are also computed by the 2D WENO type extrapolation \cite{TWSN2012}.  Then $\hat U^{(1)}(P)$ can be computed with \eqref{eq:U2}, and similarly $\hat U^{(k)}(P), k\geq 2$ can be obtained. The velocities at the ghost point can be transformed back to the original coordinate by
\begin{equation*}
\left(
\begin{array}{c}
u \\
v
\end{array}
\right)
= \mathsf{T}^{-1}
 \left(
\begin{array}{c}
\hat u  \\
\hat v
\end{array}
\right)
\end{equation*}
while the scalar quantities remain unchanged in the two coordinates.
In this way, both the solutions $U^n$ at time $t_n$ and intermediate solutions $U^{(i)}$ at ghost points are determined. Note that the 2D WENO type extrapolation has been been introduced detailedly in \cite{TWSN2012} and thus is omitted here.

%
%
%

\section{Numerical validations}

In this section we conduct several numerical experiments for both 1D and 2D problems to demonstrate the numerical stability and accuracy of our boundary treatment.
In all the computations, we use the third-order finite difference WENO scheme \cite{Liu1994weno3} with the Lax-Friedrichs flux splitting to form the numerical fluxes \cite{JS1996} and the third-order IMEX RK scheme proposed in \cite{ascher1997} with the double Butcher tableau
\begin{center}
\begin{tabular}{p{0.5cm} | p{0.8cm} p{0.8cm} p{0.8cm} p{0.8cm} p{0.8cm}}
0  &  0 & 0 & 0 & 0 & 0  \\
1/2  & 1/2 & 0 & 0 & 0 & 0  \\
2/3 & 11/18 & 1/18 & 0 & 0 & 0 \\
1/2 & 5/6 & -5/6 & 1/2 & 0 & 0 \\
1 & 1/4 & 7/4 & 3/4 & -7/4 & 0 \\
\hline
  & 1/4 & 7/4 & 3/4 & -7/4 & 0
\end{tabular}
\qquad
\begin{tabular}{p{0.5cm} | p{0.8cm} p{0.8cm} p{0.8cm} p{0.8cm} p{0.8cm}}
0   &  0  & 0 & 0 & 0 & 0  \\
1/2 &  0  & 1/2 & 0 & 0 & 0  \\
2/3 &  0  & 1/6 & 1/2 & 0 & 0 \\
1/2 &  0  & -1/2 & 1/2 & 1/2 & 0 \\
1   &  0  & 3/2 & -3/2 & 1/2 & 1/2 \\
\hline
  & 0  & 3/2 & -3/2 & 1/2 & 1/2
\end{tabular}
\end{center}
Additionally, we fix the CFL number to be $0.8$ in all the numerical examples. For all the 1D examples, the mesh for computation is taken as \eqref{3-1}, namely, the boundary is located in the middle of two grid points. Different locations of the boundary are studied for 2D problems.

\subsection{1D examples}

\noindent \textbf{Example 1.} We first consider a 1D scalar hyperbolic equation
\begin{equation}\label{528}
\partial_t u + u \partial_x u = u^2 + u,
\end{equation}
in the domain $0\le x\le 1$ with analytical solution $u=\exp(t+x)$. Since $u>0$, the boundary condition only needs to be imposed at $x=0$. Specifically, the initial and boundary conditions are set as
\begin{equation}
u(x,0) = \exp(x), \qquad
u(0,t) = \exp(t).
\end{equation}
We compute the $L^1, L^2$ and $L^{\infty}$ errors at the terminal time $t=1$, and the results are listed in Table \ref{tb:example1}. It can be seen that the convergence order is around three.

\begin{table}[!htbp]\centering
\caption{Example 1: error table for the scalar nonlinear equation \eqref{528}. }
\begin{tabular}{lllllll}
\toprule
$\Delta x$ & $L^1$ error & order  & $L^2$ error & order  & $L^{\infty}$ error   & order\\ \hline
1/20      &  7.62e-4    &         & 1.29e-3     &        & 4.22e-3   &      \\
1/40      &  1.13e-4    & 2.75    & 2.17e-4     &  2.57  & 9.57e-4   &  2.14\\
1/80      &  1.41e-5    & 3.01    & 2.59e-5     &  3.07  & 1.07e-4   &  3.16\\
1/160     &  1.50e-6    & 3.23    & 2.42e-6     &  3.42  & 7.20e-6   &  3.89\\
1/320     &  1.57e-7    & 3.25    & 2.32e-7     &  3.38  & 5.53e-7   &  3.70\\
\bottomrule
\end{tabular}\label{tb:example1}
\end{table}

\noindent \textbf{Example 2.} Next we solve a 1D linear system
\begin{equation} \label{415}
\begin{split}
& \partial_t u + \partial_x v = 0 \\
& \partial_t v + \partial_x u = -\frac{1}{\epsilon}(u+v)
\end{split}
\end{equation}
in the domain $0\le x\le 1$. The Jacobian matrix has two eigenvalues $\pm1$ and thus only one boundary condition should be imposed at $x=0$ and also $x=1$, respectively.
We construct two sets of analytical solutions to \eqref{415}:
\begin{equation} \label{416}
u(t,x) =\exp(t+x), \quad  v(t,x) = -\exp(t+x)
\end{equation}
and
\begin{equation} \label{417}
u(t,x) =\exp(t+x) + \exp(-\frac{x}{\epsilon} ), \qquad  v(t,x) = -\exp(t+x).
\end{equation}
There exists no boundary layer for the first solution \eqref{416} while the boundary layer occurs at $x=0$ for the second one \eqref{417} with small $\epsilon>0$. We test our method for the two cases by taking the initial data to match the analytical solution. The mesh is given as \eqref{25} and the solutions of $u$ at $t=1$ are used to evaluate the errors. In the second case, the three interior grid points near the left boundary are excluded in the computation of the errors due to the existence of the boundary layer.

\begin{figure}[!ht]
\centering
\includegraphics[width=0.45\textwidth]{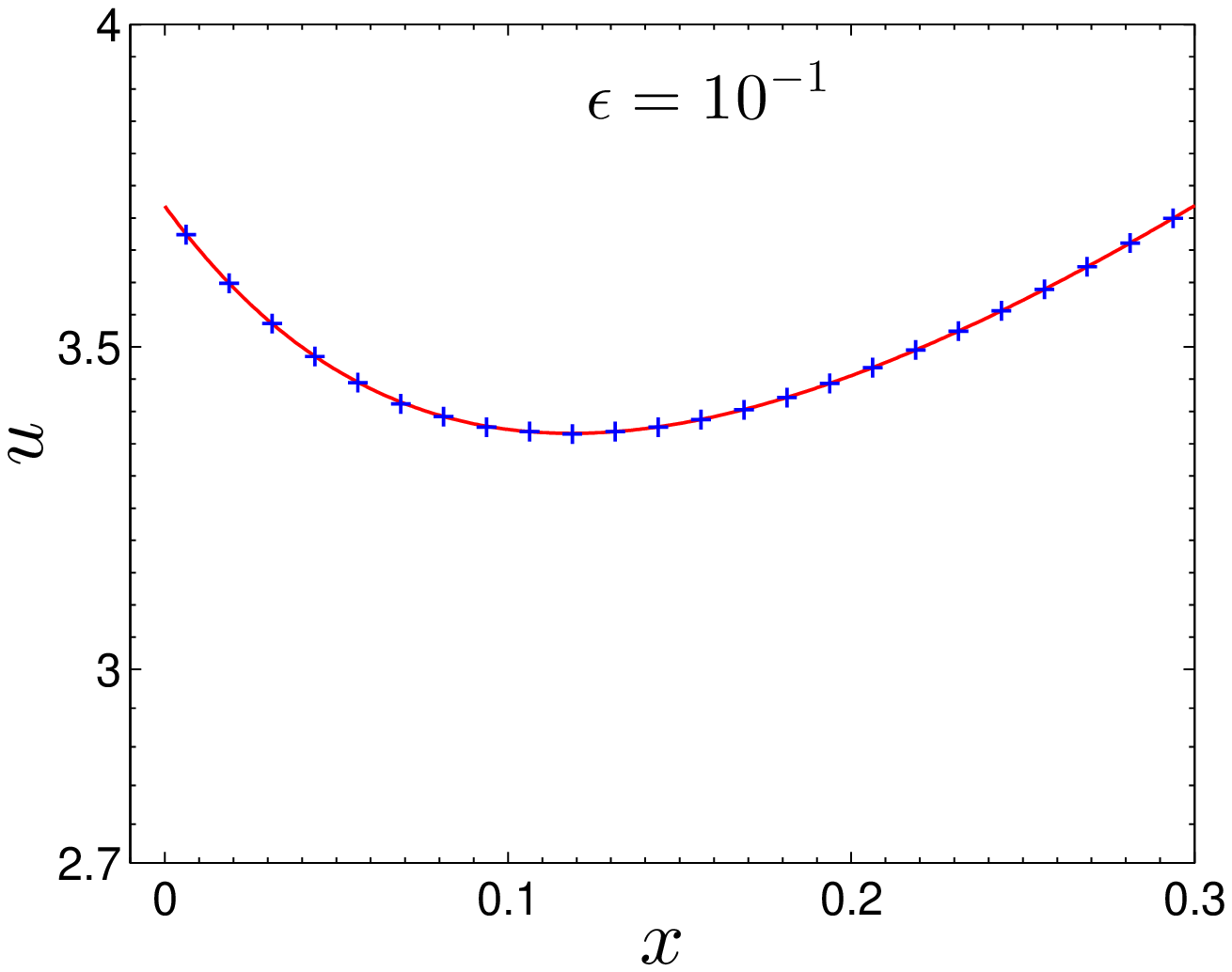}
\includegraphics[width=0.45\textwidth]{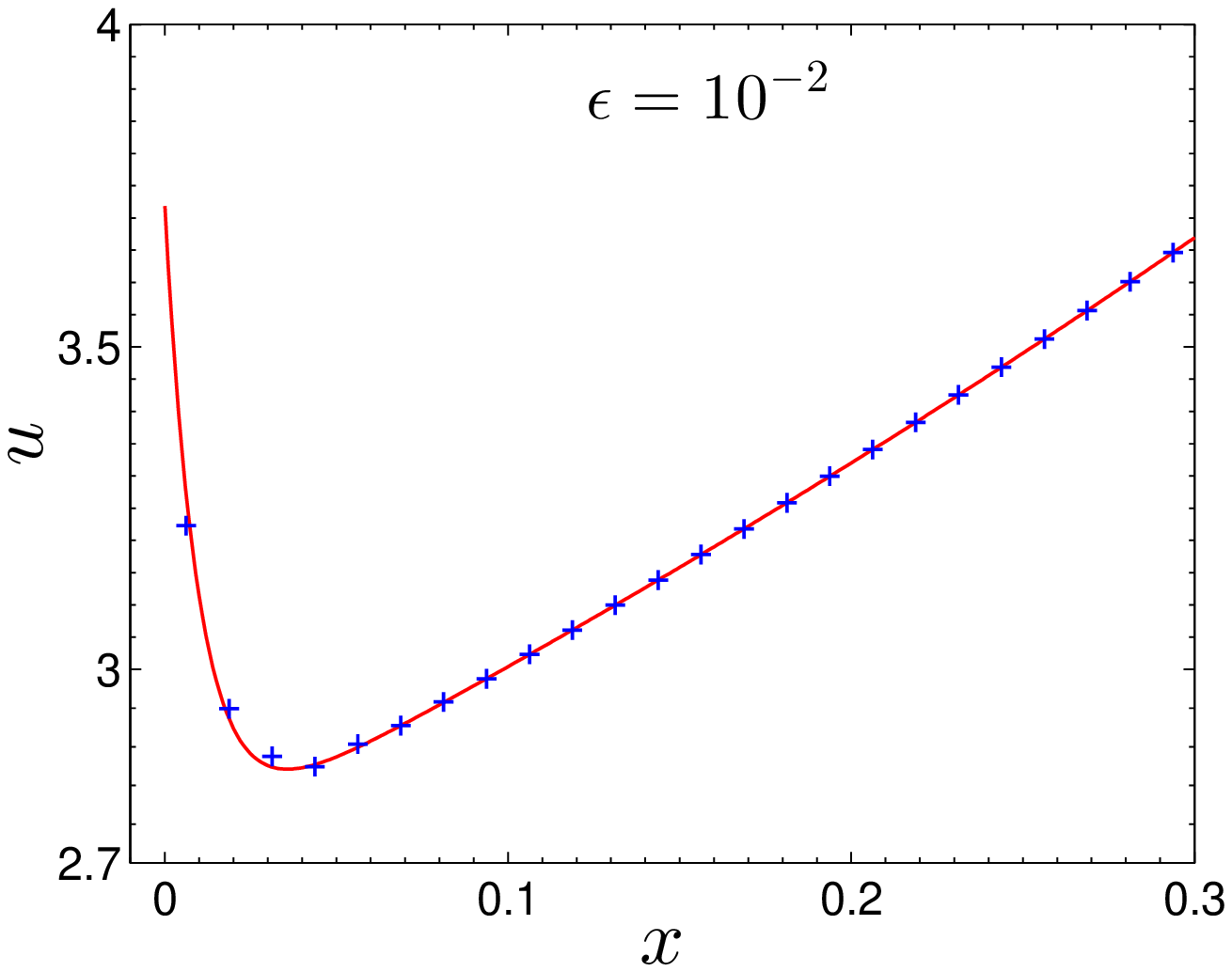}
\includegraphics[width=0.45\textwidth]{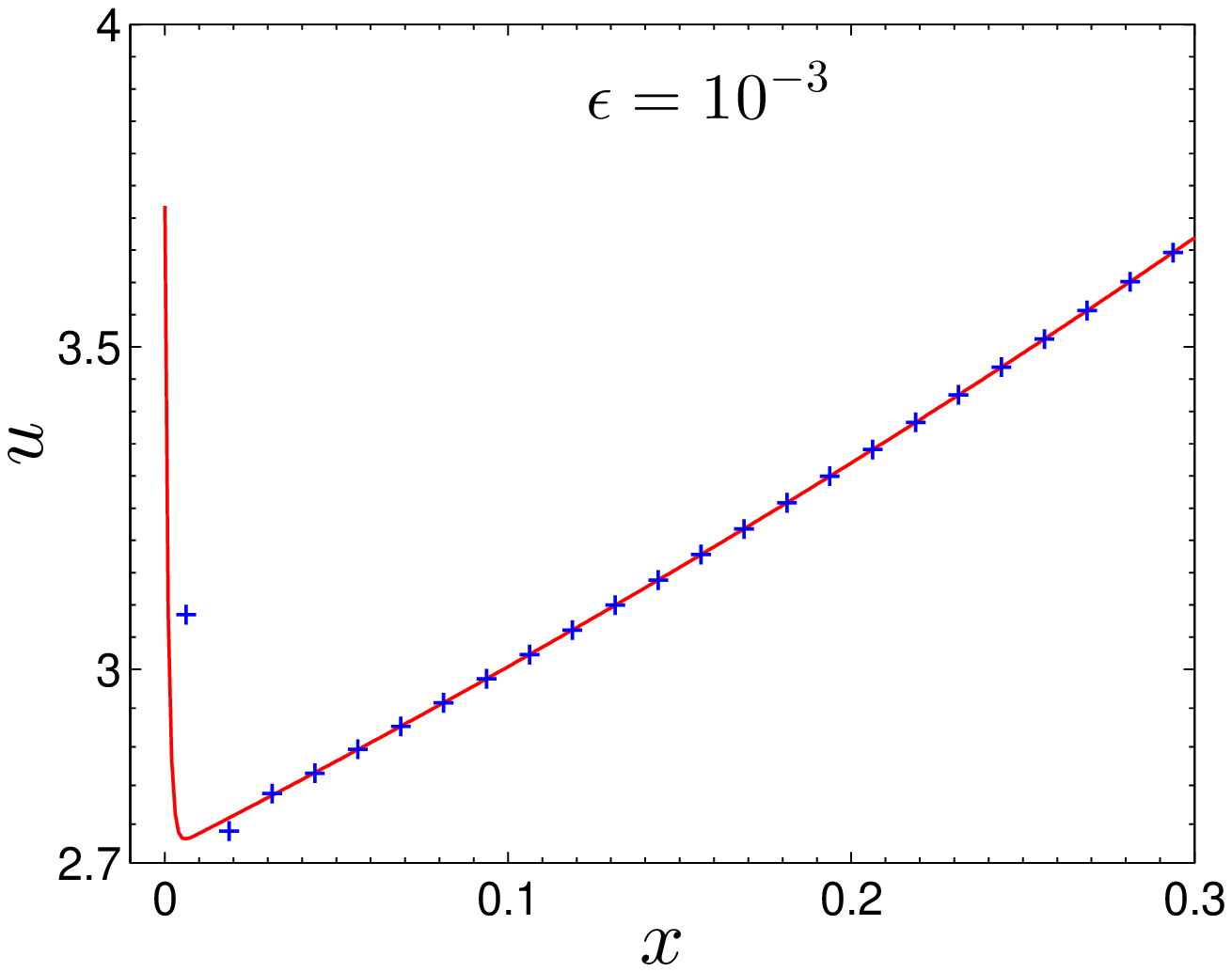}
\includegraphics[width=0.45\textwidth]{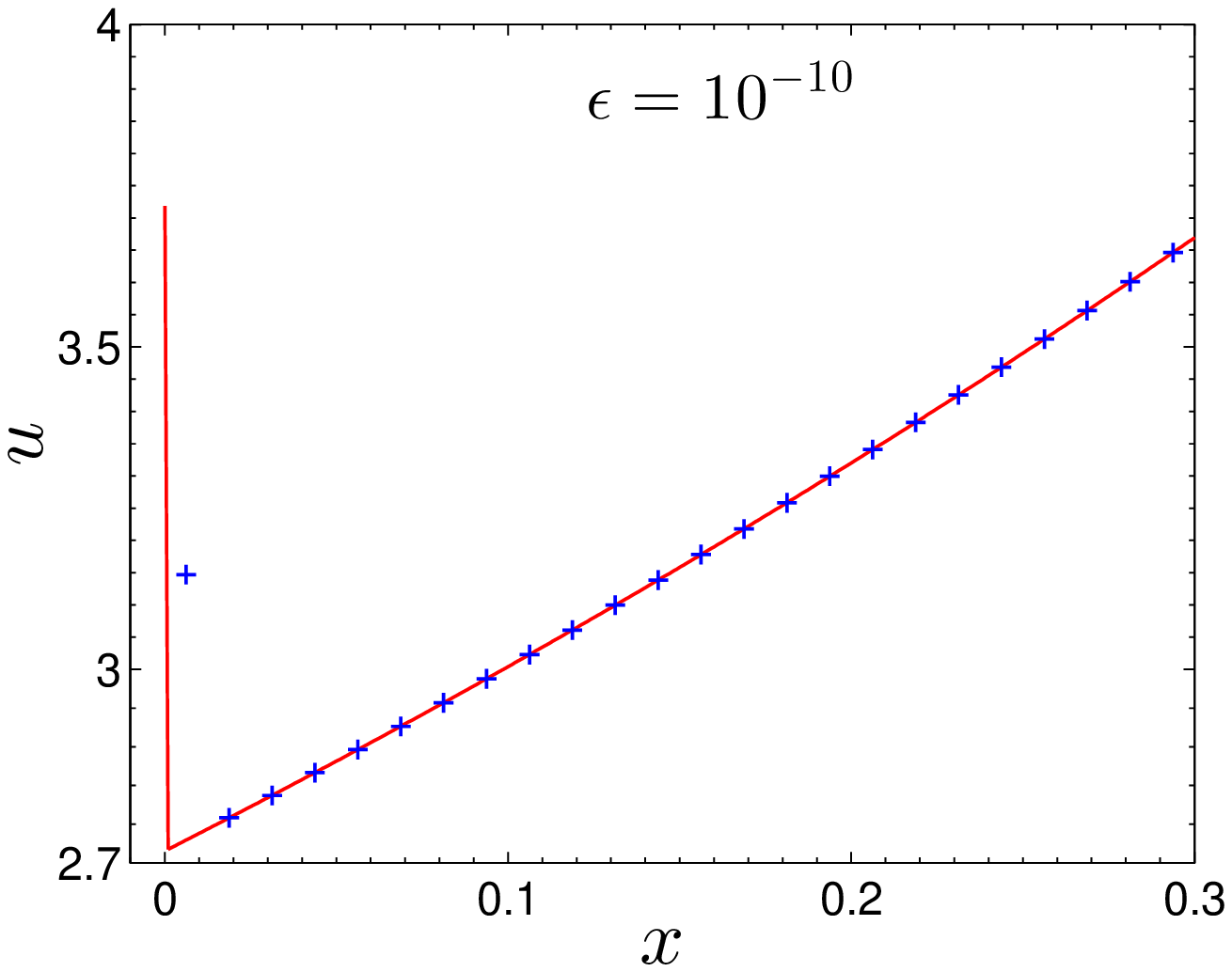}
\caption{Example 2: Comparison of the computational results (markers) and analytical solution (line) of the linear system \eqref{415} with different $\epsilon$ and fixed $\Delta x = 1/80$.}
\label{fig:1D_example2-1}
\end{figure}

For the solution \eqref{416}, we specify $u$ and $v$ at the left and right boundaries, respectively, \ie,
\begin{equation} \label{418_0}
\begin{split}
& u(t,0) = \exp(t) \quad \mbox{at} \quad x=0, \\
& v(t,1) = -\exp(t+1)\quad \mbox{at} \quad x=1.
\end{split}
\end{equation}
Similarly, the boundary condition for \eqref{417} is given by
\begin{equation} \label{418}
\begin{split}
& u(t,0) = \exp(t)+ 1 \quad \mbox{at} \quad x=0, \\
& v(t,1) = -\exp(t+1)\quad \mbox{at} \quad x=1.
\end{split}
\end{equation}
The errors for the first case with solution \eqref{416} are given in Table \ref{tb:example2-1}, from which the desired third-order convergence can be observed.
For the second case with solution \eqref{417}, we take $\epsilon=1$ and $10^{-10}$ to test the convergence order. The results in Table \ref{tb:example2-2} show that the present method of boundary treatment is third-order accurate for both choices of $\Delta x\ll \epsilon$ and $\Delta x\gg \epsilon$.
Fig.~\ref{fig:1D_example2-1} shows the solutions near the left boundary for different $\epsilon$ with $\Delta x = 1/80$. It can be seen that our method behaves well even in the presence of boundary layer.


\begin{table}[!htbp]\centering
\caption{Example 2: error table for the linear system \eqref{415} with solution \eqref{416}. }
\begin{tabular}{lllllll}
\toprule
$\Delta x$ & $L^1$ error & order  & $L^2$ error & order  & $L^{\infty}$ error   & order\\ \hline
1/20      &  9.77e-4    &         & 1.06e-3     &        & 2.37e-3   &      \\
1/40      &  1.08e-4    & 3.03    & 1.25e-4     &  3.08  & 3.20e-4   &  2.89\\
1/80      &  1.46e-5    & 2.88    & 1.71e-5     &  2.87  & 6.59e-5   &  2.28\\
1/160     &  1.78e-6    & 3.03    & 1.95e-6     &  3.14  & 5.58e-6   &  3.56\\
1/320     &  2.35e-7    & 2.93    & 2.52e-7     &  2.95  & 7.24e-7   &  2.95\\
\bottomrule
\end{tabular}\label{tb:example2-1}
\end{table}

\begin{table}[!htbp]\centering
\caption{Example 2: Error table for the linear system \eqref{415} with solution \eqref{417}. }
\begin{tabular}{lllllll}
\toprule
         &             &         &  $\epsilon = 1$  &           &  \\\hline
$\Delta x$ & $L^1$ error & order  & $L^2$ error & order  & $L^{\infty}$ error   & order\\ \hline
1/20      &  1.01e-3    &         & 1.18e-3     &        & 2.39e-3   &      \\
1/40      &  1.20e-4    & 3.08    & 1.34e-4     &  3.14  & 3.17e-4   &  2.91\\
1/80      &  1.57e-5    & 2.93    & 1.78e-5     &  2.91  & 6.49e-5   &  2.29\\
1/160     &  1.90e-6    & 3.05    & 2.02e-6     &  3.14  & 5.52e-6   &  3.56\\
1/320     &  2.49e-7    & 2.93    & 2.62e-7     &  2.95  & 7.19e-7   &  2.94\\
\toprule
          &             &         &  $\epsilon = 10^{-10}$  &           &  \\\hline
$\Delta x$ & $L^1$ error & order  & $L^2$ error & order  & $L^{\infty}$ error   & order\\ \hline
1/20      &  9.66e-4    &         & 1.16e-3     &        & 2.38e-3   &      \\
1/40      &  1.16e-4    & 3.06    & 1.32e-4     &  3.13  & 3.22e-4   &  2.89\\
1/80      &  1.53e-5    & 2.93    & 1.77e-5     &  2.90  & 6.58e-5   &  2.29\\
1/160     &  1.84e-6    & 3.04    & 2.01e-6     &  3.14  & 5.72e-6   &  3.52\\
1/320     &  2.51e-7    & 2.88    & 2.68e-7     &  2.91  & 7.07e-7   &  3.02\\
\bottomrule
\end{tabular}\label{tb:example2-2}
\end{table}

\noindent \textbf{Example 3.} We further construct the following nonlinear system
\begin{equation} \label{420}
\begin{split}
& \partial_t u + \partial_x v = 0 \\
& \partial_t v + \partial_x [u+\frac{1}{2}(u+v)^2] = -\frac{1}{\epsilon}[u+v+(u+v)^2]
\end{split}
\end{equation}
in the domain $0\le x\le1$ to test the accuracy of our method. The expressions \eqref{416} and \eqref{417} are also analytical solutions to \eqref{420} if cooperated with appropriate initial and boundary conditions. The Jacobian matrix of the flux is
\begin{equation*}
A =
\left(
\begin{array}{cc}
 0 & 1\\
 1+u+v & u+v
\end{array}
\right)
\end{equation*}
and it can be diagonalized as
\begin{equation*}
LAR =
\left(
\begin{array}{cc}
 1+u+v & 0\\
0 & -1
\end{array}
\right) ,
\end{equation*}
where $1+u+v>0$ and
\begin{equation*}
L=
\left(
\begin{array}{cc}
1 & 1 \\
1 & -\frac{1}{1+u+v}
\end{array}
\right), \quad
R=L^{-1}=
\left(
\begin{array}{cc}
\frac{1}{2+u+v} & \frac{1+u+v}{2+u+v} \\
\frac{1+u+v}{2+u+v} & -\frac{1+u+v}{2+u+v}
\end{array}
\right).
\end{equation*}
The boundary conditions in \eqref{418_0} and \eqref{418} are employed here.

The errors of the solution \eqref{416} with boundary condition \eqref{418_0} are given in Table \ref{tb:example3-1} and the desired third-order accuracy is achieved again. Similar to the linear system, we take $\epsilon=1$ and $10^{-10}$ for the  solution \eqref{417} with boundary condition \eqref{418} and the errors are shown in Table \ref{tb:example3-2}. It is clear that the convergence orders are about 3 for both choices of $\epsilon$. These results demonstrate the third-order accuracy of our method for nonlinear systems.

\begin{table}[!htbp]\centering
\caption{Example 3: error table for the nonlinear system \eqref{420} with solution \eqref{416}. }
\begin{tabular}{lllllll}
\toprule
$\Delta x$ & $L^1$ error & order  & $L^2$ error & order  & $L^{\infty}$ error   & order\\ \hline
1/20      &  9.25e-4    &         & 1.09e-3     &        & 2.39e-3   &      \\
1/40      &  1.13e-4    & 3.03    & 1.29e-4     &  3.08  & 3.22e-4   &  2.89\\
1/80      &  1.52e-5    & 2.89    & 1.76e-5     &  2.87  & 6.61e-5   &  2.28\\
1/160     &  1.86e-6    & 3.03    & 2.01e-6     &  3.13  & 5.61e-6   &  3.56\\
1/320     &  2.44e-7    & 2.93    & 2.60e-7     &  2.95  & 7.28e-7   &  2.95\\
\bottomrule
\end{tabular}\label{tb:example3-1}
\end{table}

\begin{table}[!htbp]\centering
\caption{Example 3: error table for the nonlinear system \eqref{420} with solution \eqref{417}. }
\begin{tabular}{lllllll}
\toprule
          &             &         &  $\epsilon = 1$  &           &  \\\hline
$\Delta x$ & $L^1$ error & order  & $L^2$ error & order  & $L^{\infty}$ error   & order\\ \hline
1/20      &  2.52e-4    &         & 2.99e-4     &        & 4.33e-4   &      \\
1/40      &  4.28e-5    & 2.56    & 4.77e-5     &  2.65  & 6.74e-5   &  2.68\\
1/80      &  6.09e-6    & 2.81    & 6.58e-6     &  2.86  & 8.98e-6   &  2.91\\
1/160     &  7.86e-7    & 2.95    & 8.33e-7     &  2.98  & 1.10e-6   &  3.03\\
1/320     &  9.25e-8    & 3.09    & 9.66e-8     &  3.11  & 1.22e-7   &  3.18\\
\toprule
          &             &         &  $\epsilon = 10^{-10}$  &           &  \\\hline
$\Delta x$ & $L^1$ error & order  & $L^2$ error & order  & $L^{\infty}$ error   & order\\ \hline
1/20      &  9.66e-4    &         & 1.16e-3     &        & 2.38e-3   &      \\
1/40      &  1.16e-4    & 3.06    & 1.32e-4     &  3.13  & 3.22e-4   &  2.89\\
1/80      &  1.53e-5    & 2.93    & 1.77e-5     &  2.90  & 6.58e-5   &  2.29\\
1/160     &  1.84e-6    & 3.04    & 2.01e-6     &  3.14  & 5.72e-6   &  3.52\\
1/320     &  2.51e-7    & 2.88    & 2.68e-7     &  2.91  & 7.07e-7   &  3.02\\
\bottomrule
\end{tabular}\label{tb:example3-2}
\end{table}

\subsection{2D reactive Euler equations}

 All the following examples are devoted to validate our boundary treatment for the 2D reactive Euler equations \eqref{eq:euler}. The parameters are fixed as in \cite{Wang2012JCP}:
\[
\gamma = 1.2, \quad q=50, \quad \tilde T = 50, \quad \tilde K=2566.4.
\]

\begin{figure}[!ht]
\centering
\includegraphics[width=0.45\textwidth]{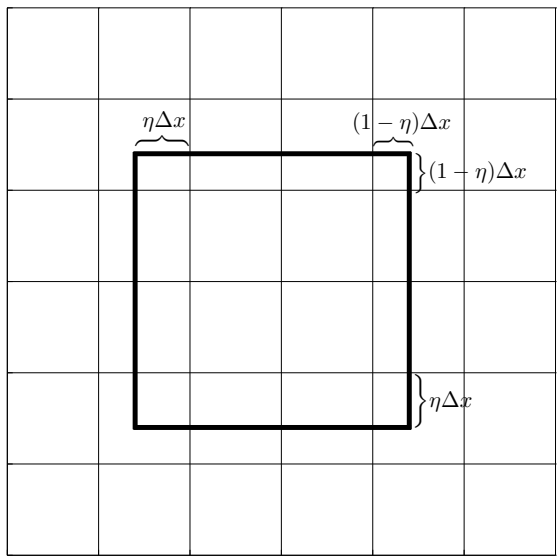}
\caption{Example 4: Illustration of the mesh. The thin lines are the grid lines and the thick lines are the boundaries of the square domain $[0,1]\times[0,1]$. The distance between the bottom boundary and its adjacent gird line is $\eta \Delta x$ and that at the top boundary is $(1-\eta) \Delta x$. In a similar way the left and right boundaries are placed. Two layers of grid points are required outside each boundary for the third-order WENO scheme.}
\label{fig:2D_example0_mesh}
\end{figure}

\noindent \textbf{Example 4.}
We first verify the accuracy of the method through a problem with analytical solution
\begin{equation}\label{2D_example0_solution}
\rho = 1+0.3\sin[ 2\pi(x+y-t) ], \quad u=1, \quad v=0, \quad p=1, \quad Y=0
\end{equation}
in a square domain $[0,1]\times[0,1]$, subject to periodic boundary conditions at the left and right boundaries and $v=0$ at the top and bottom boundaries. The latter boundary condition is treated with our method. The mesh for computation is illustrated in Fig.~\ref{fig:2D_example0_mesh}. As shown in the figure, the distance between the bottom boundary and its adjacent gird line is $\eta \Delta x$ and that at the top boundary is $(1-\eta) \Delta x$. The left and right boundaries are placed in a similar way. We test the accuracy with $\eta = 0.5, 0.7$ and 0.9, and the results are listed in Table \ref{tb:2D_example0}. It can be seen that the influence of $\eta$ on the errors is slight and the convergence order is about 3 for all the three choices of $\eta$. These demonstrate the third-order accuracy of our method  for 2D problems.

\begin{table}[!htbp]\centering
\caption{Example 4: error table for 2D reactive euler equations with solution \eqref{2D_example0_solution}.}
\begin{tabular}{lllllll}
\toprule
          &             &         &  $\eta = 0.5$  &           &  \\\hline
$\Delta x$ & $L^1$ error & order  & $L^2$ error & order  & $L^{\infty}$ error   & order\\ \hline
1/20      &  1.37e-2    &         & 1.64e-2     &        & 3.18e-2   &      \\
1/40      &  4.12e-3    & 1.73    & 5.55e-3     &  1.56  & 1.29e-2   &  1.30\\
1/80      &  7.78e-4    & 2.40    & 1.27e-3     &  2.13  & 3.75e-3   &  1.78\\
1/160     &  7.35e-5    & 3.40    & 1.35e-4     &  3.23  & 5.25e-4   &  2.84\\
\toprule
          &             &         &  $\eta = 0.7$  &           &  \\\hline
$\Delta x$ & $L^1$ error & order  & $L^2$ error & order  & $L^{\infty}$ error   & order\\ \hline
1/20      &  1.41e-2    &         & 1.64e-2     &        & 3.10e-2   &      \\
1/40      &  4.10e-3    & 1.78    & 5.56e-3     &  1.56  & 1.27e-2   &  1.28\\
1/80      &  7.79e-4    & 2.40    & 1.27e-3     &  2.13  & 3.76e-3   &  1.76\\
1/160     &  7.34e-5    & 3.41    & 1.35e-4     &  3.23  & 5.32e-4   &  2.82\\
\toprule
          &             &         &  $\eta = 0.9$  &           &  \\\hline
$\Delta x$ & $L^1$ error & order  & $L^2$ error & order  & $L^{\infty}$ error   & order\\ \hline
1/20      &  1.44e-2    &         & 1.66e-2     &        & 3.08e-2   &      \\
1/40      &  4.12e-3    & 1.81    & 5.58e-3     &  1.57  & 1.26e-2   &  1.28\\
1/80      &  7.70e-4    & 2.42    & 1.27e-3     &  2.14  & 3.69e-3   &  1.78\\
1/160     &  7.34e-5    & 3.39    & 1.35e-4     &  3.23  & 5.27e-4   &  2.81\\
\bottomrule
\end{tabular}\label{tb:2D_example0}
\end{table}

\noindent \textbf{Example 5.}
Next we test our method with the example in \cite{Wang2012JCP} containing discontinuity. The computational domain is $[0,2]\times [0,2]$ and the initial data is
\begin{equation}\label{2D_example0_InitialCondition}
(\rho, u, v, p, Y)
=
\left\{
\begin{array}{l}
( 1,0,0,80, 0.2 ), \quad x^2+y^2 \leq 0.36,  \\
( 1,0,0,10, 0.8 ), \quad \mbox{otherwise},
\end{array}
\right.
\end{equation}
where $p$ and $Y$ are different from that in \cite{Wang2012JCP} to avoid possibly negative pressure or density.
The boundary condition for the left and bottom boundaries are $u=0$ and $v=0$, respectively, and the solutions of the grid points outside the right and top boundaries are equal to those at the grid points inside the domain.

\begin{figure}[!ht]
\centering
\includegraphics[width=0.45\textwidth]{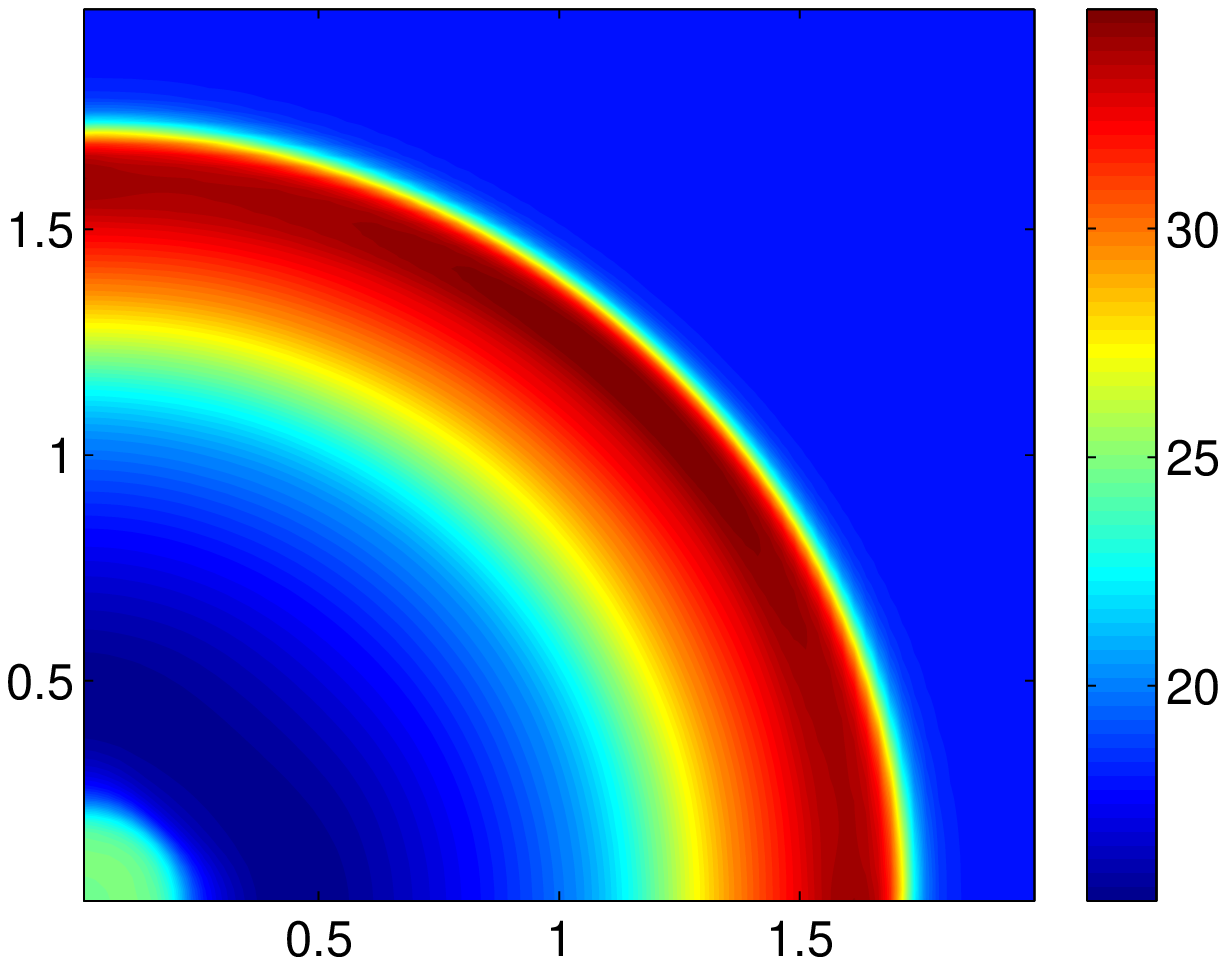}
\includegraphics[width=0.45\textwidth]{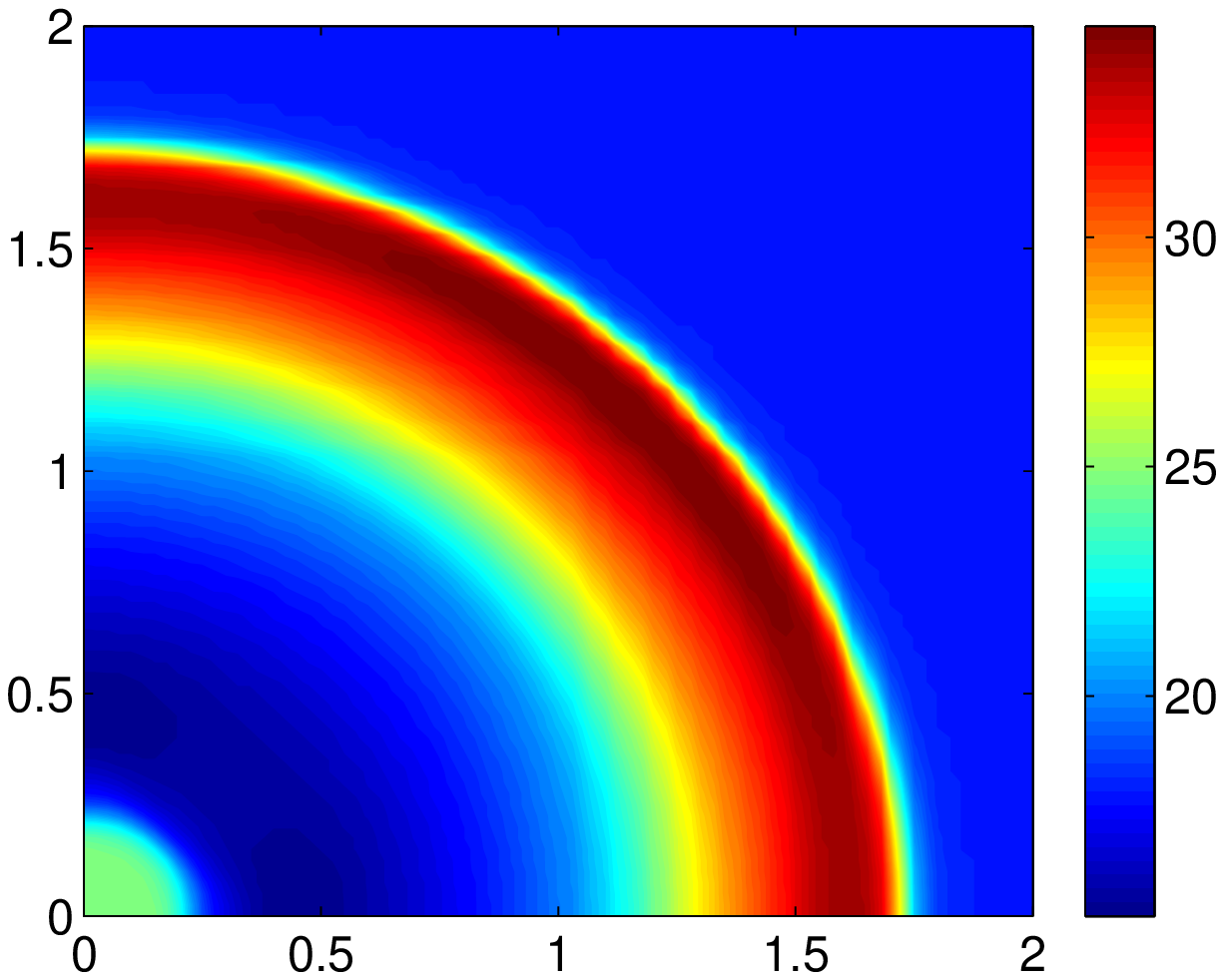}
\includegraphics[width=0.45\textwidth]{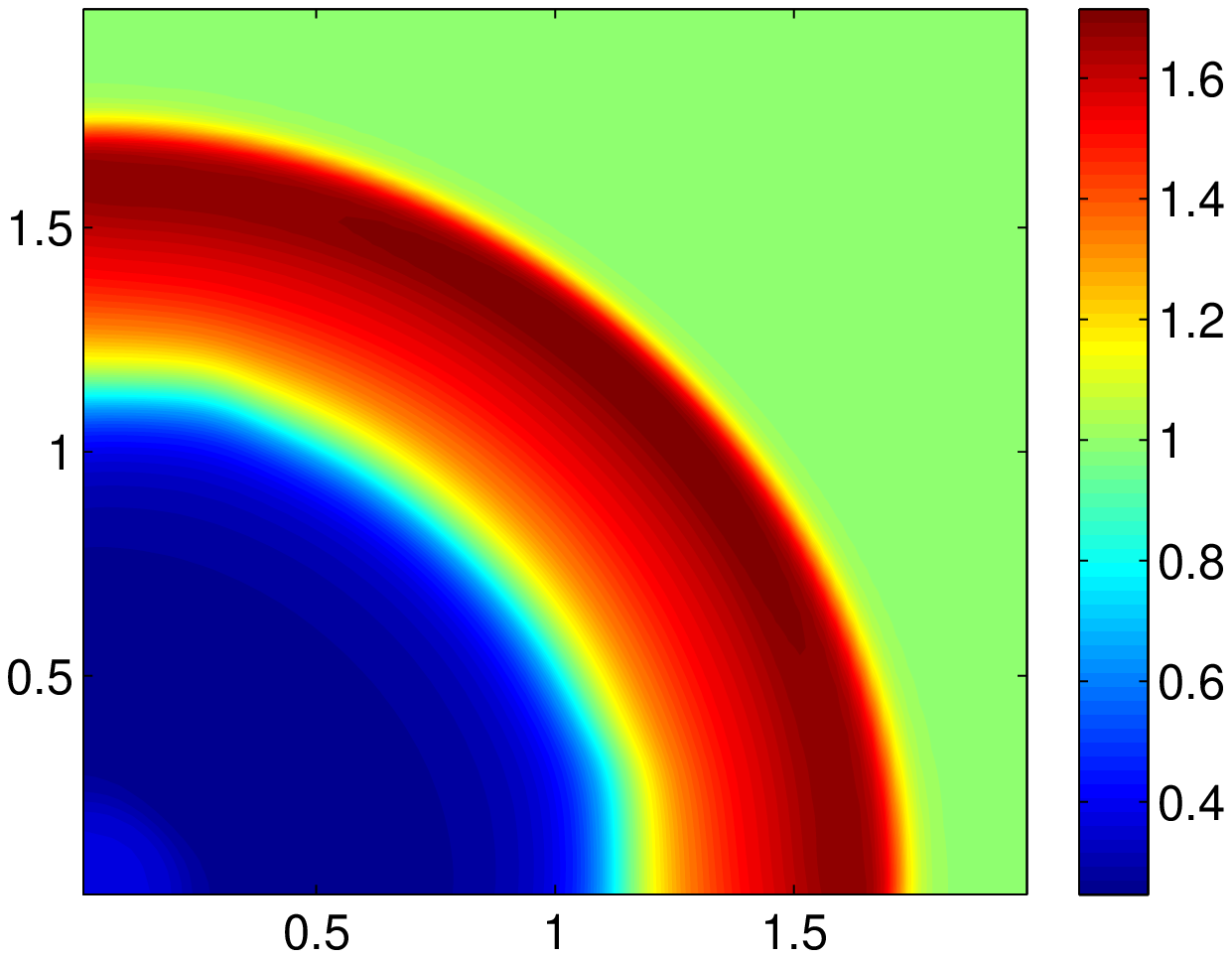}
\includegraphics[width=0.45\textwidth]{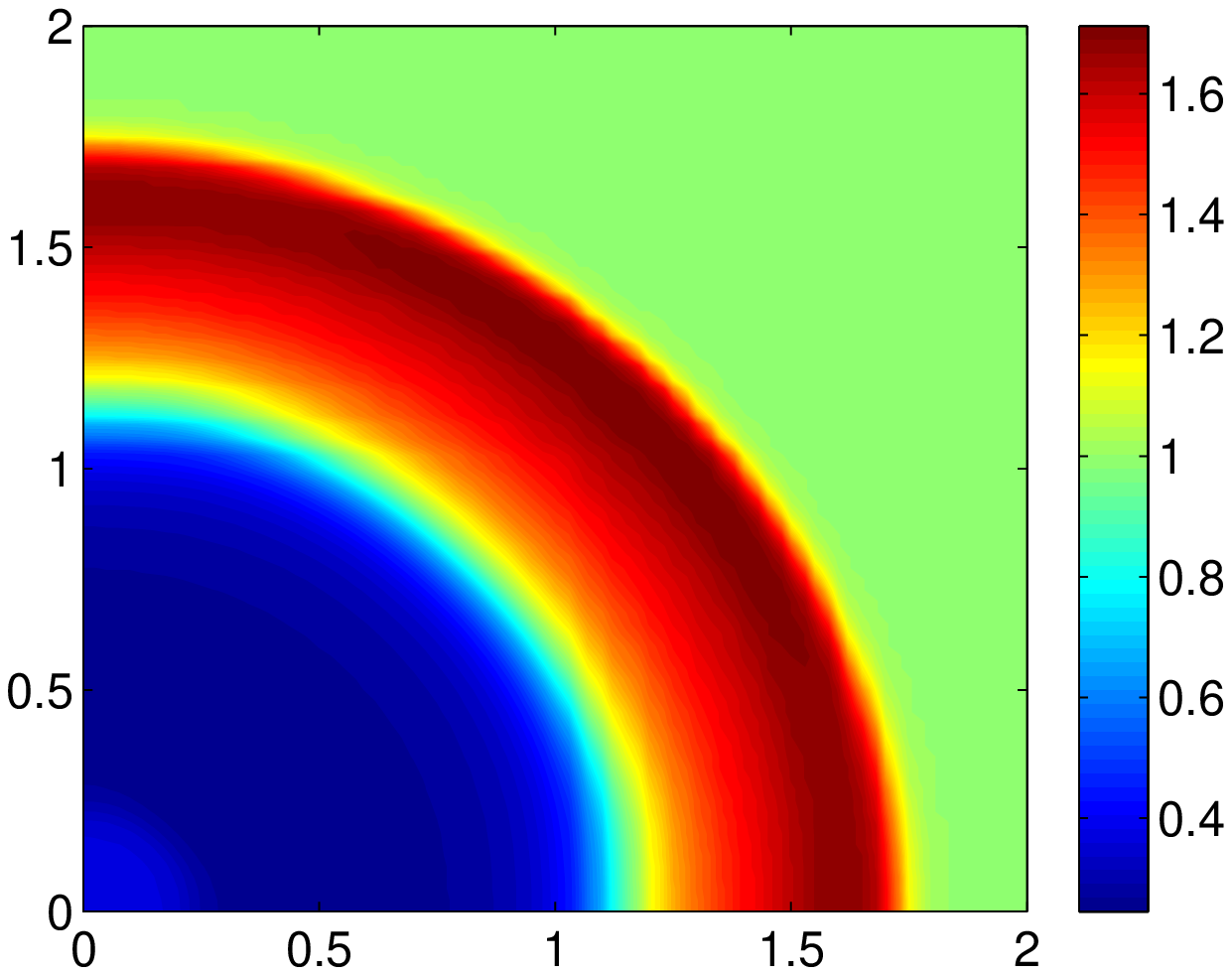}
\caption{Example 5: Contour plots of pressure (top) and density (bottom) of the present method (left) and the third-order positivity-preserving scheme in \cite{Huang2019jsc} (right) with $\Delta x = 1/40$ at $t=0.16$.}
\label{fig:2D_example1_p_rho_contour}
\end{figure}

\begin{figure}[!ht]
\centering
\includegraphics[width=0.45\textwidth]{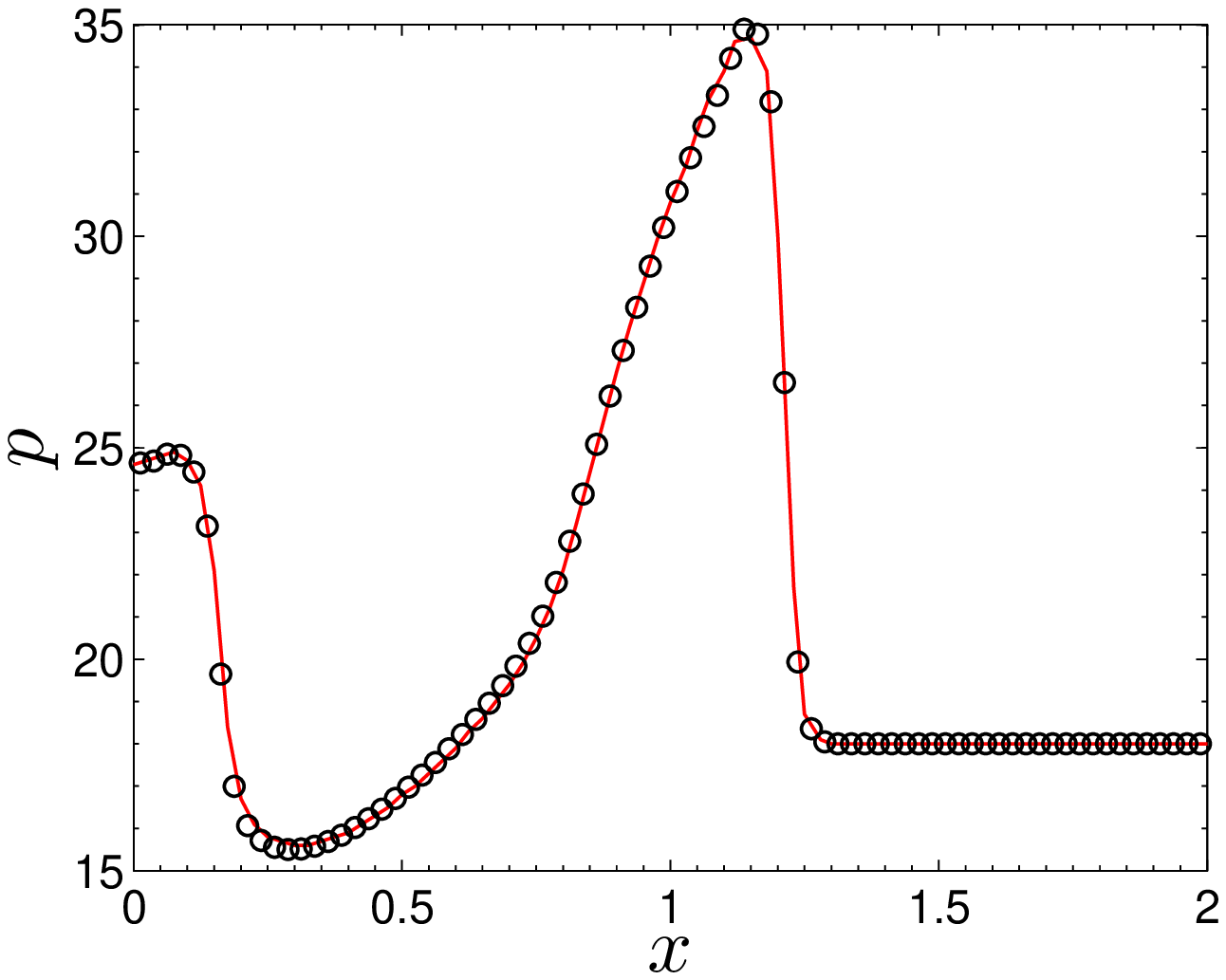}
\includegraphics[width=0.45\textwidth]{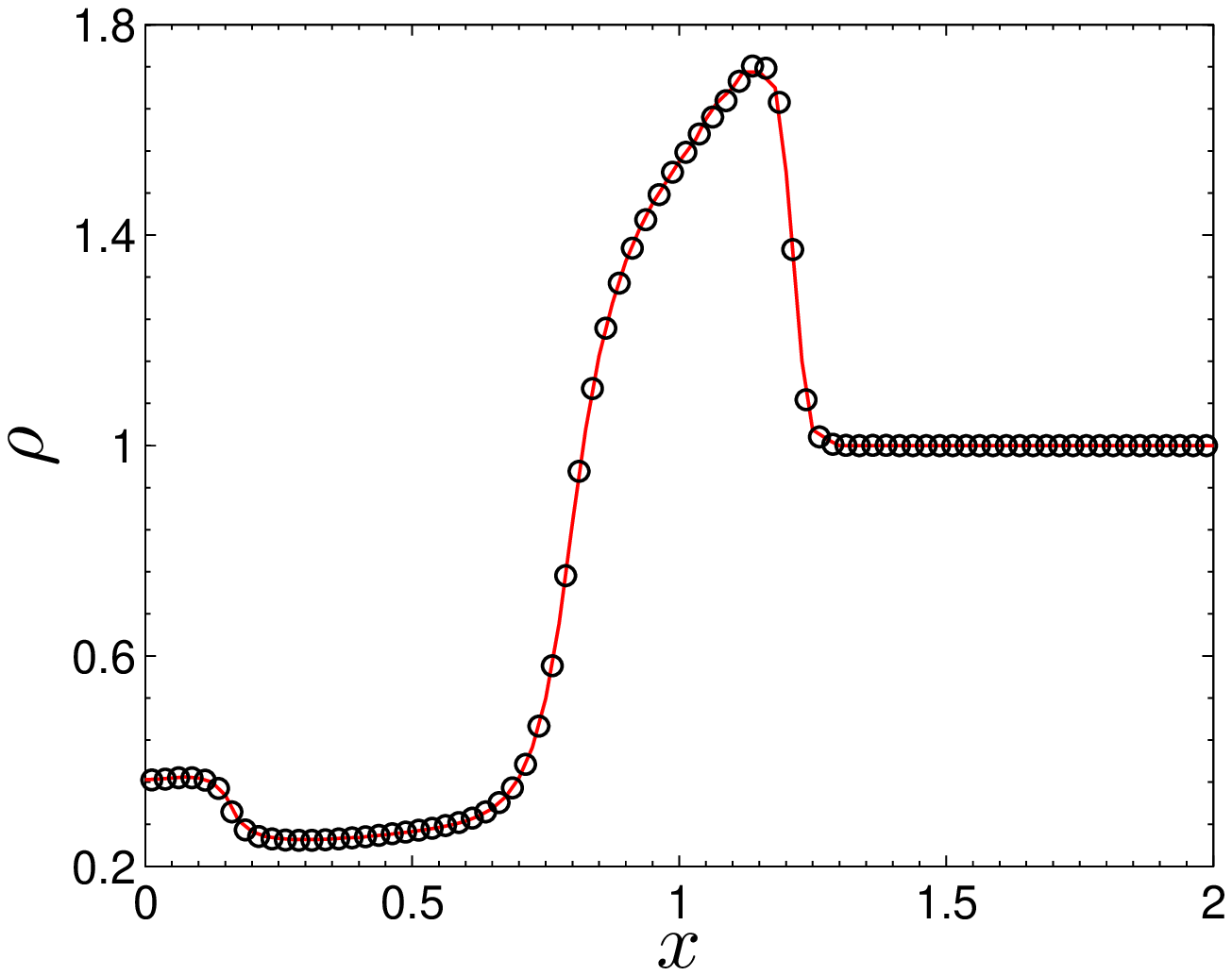}
\caption{Example 5: Comparison of the pressure (left) and density (right) of the present method (circle) and the third-order positivity-preserving scheme in \cite{Huang2019jsc} (line) along the diagonal line $y=x$.}
\label{fig:2D_example1_p_rho_diagonal}
\end{figure}

\begin{figure}[!ht]
\centering
\includegraphics[width=0.45\textwidth]{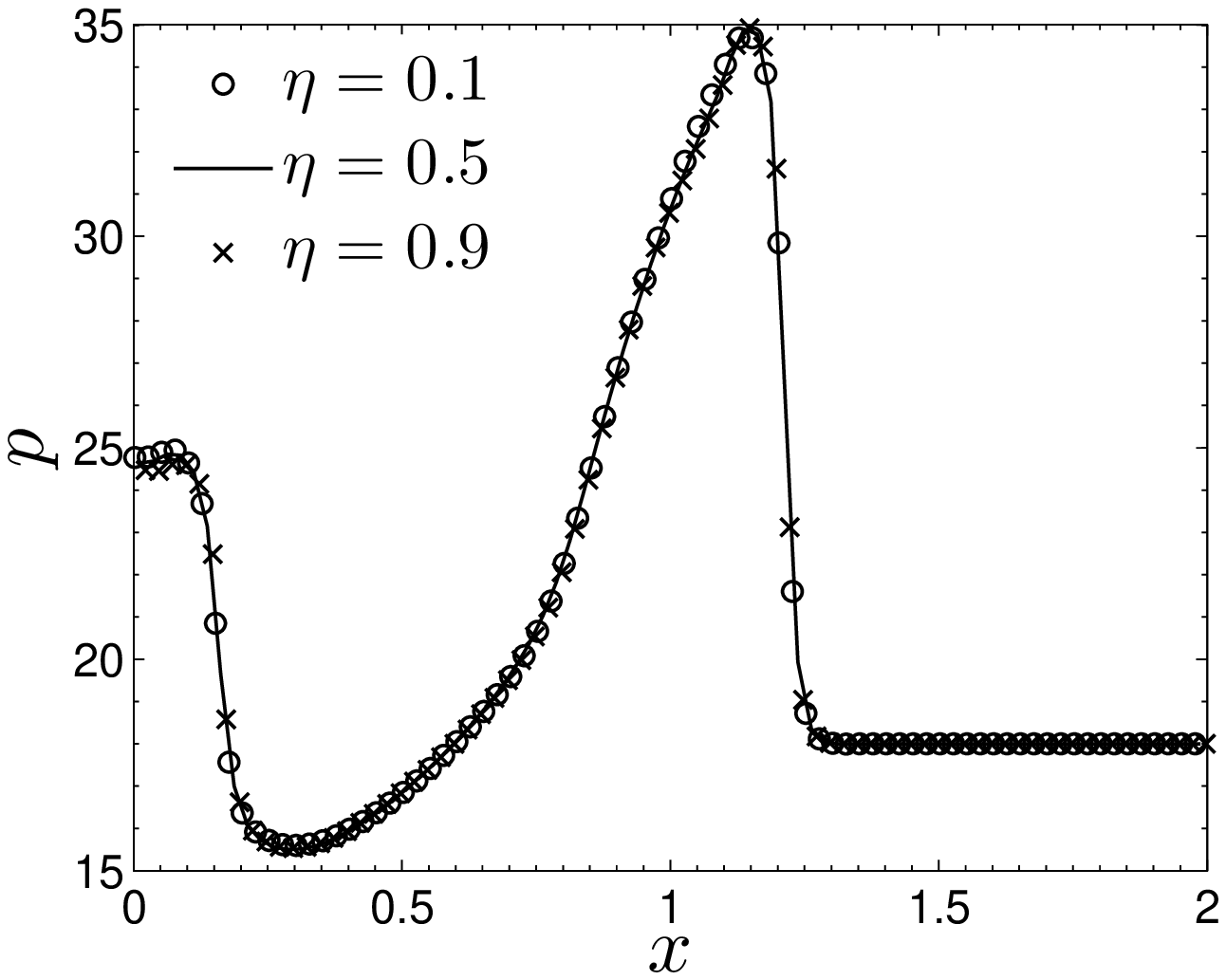}
\includegraphics[width=0.45\textwidth]{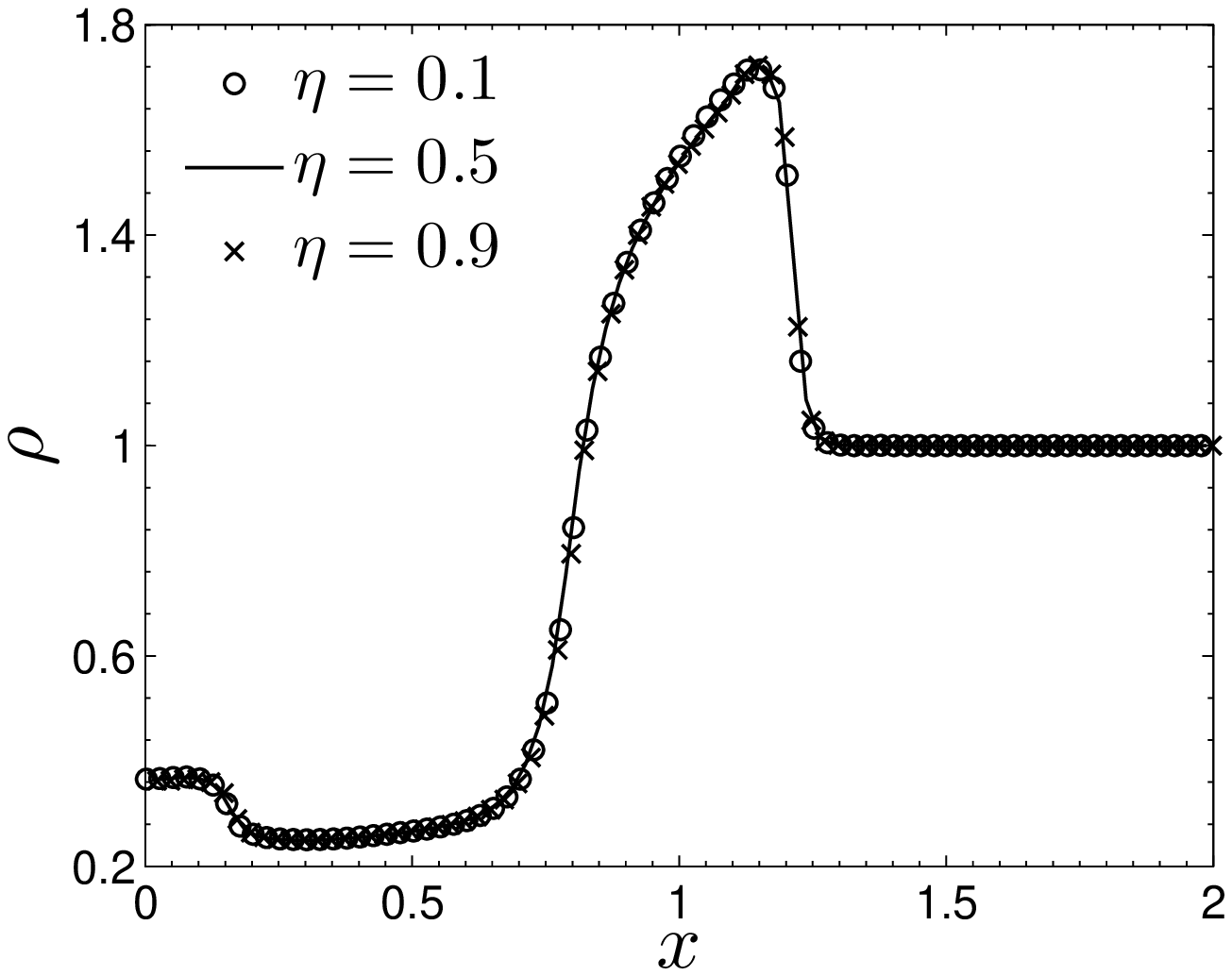}
\caption{Example 5: Comparison of the pressure (left) and density (right) of the present method with different $\eta$ along the diagonal line $y=x$. }
\label{fig:2D_example1_p_rho_DifferentEta}
\end{figure}

The mesh for computation is the same as that illustrated in Fig.~\ref{fig:2D_example0_mesh} with a ratio $\eta$.
We compare our results with those of the third-order positivity-preserving scheme in \cite{Huang2019jsc} with reflective boundary conditions (left and bottom) under the same mesh size $\Delta x = 1/40$.
The contour plots of $p$ and $\rho$ at $t=0.16$ are given in Fig.~\ref{fig:2D_example1_p_rho_contour} and the solutions along the horizontal diagonal line $y=x$ are plotted in Fig.~\ref{fig:2D_example1_p_rho_diagonal}. Good agreement between the two methods can be observed.
Additionally, we compare in Fig.~\ref{fig:2D_example1_p_rho_DifferentEta} the solutions of different ratio $\eta$.  It can be seen again that the influence of $\eta$ is small.

\noindent \textbf{Example 6.}
We test the detonation diffraction in this example with the same configuration of boundaries in \cite{Wang2012JCP} (see Fig.~\ref{fig:2D_example2_p_rho_contour} below).
The spatial domain is $[0,5] \times [0,5] $ and the initial condition is
\begin{equation*}
(\rho, u, v, E, Y)
=
\left\{
\begin{array}{l}
( 11,6.18,0,970, 1 ), \quad x<0.5,  \\
( 5,0,0,400, 1), \quad \mbox{otherwise},
\end{array}
\right.
\end{equation*}
which is constructed based on that in \cite{Wang2012JCP} to avoid negative pressure or density.
The boundary condition is $\bm n \cdot \bm u = 0$ with $\bm n$ being the normal direction of the boundary, except that at $x=0$, $(\rho, u, v, E, Y) = ( 11,6.18,0,970, 1 )$. The terminal time is $t=0.6$ and the mesh size is $\Delta x = 1/20$. The ratio $\eta$ is fixed as 0.5 for all the boundaries.

\begin{figure}[!ht]
\centering
\includegraphics[width=0.45\textwidth]{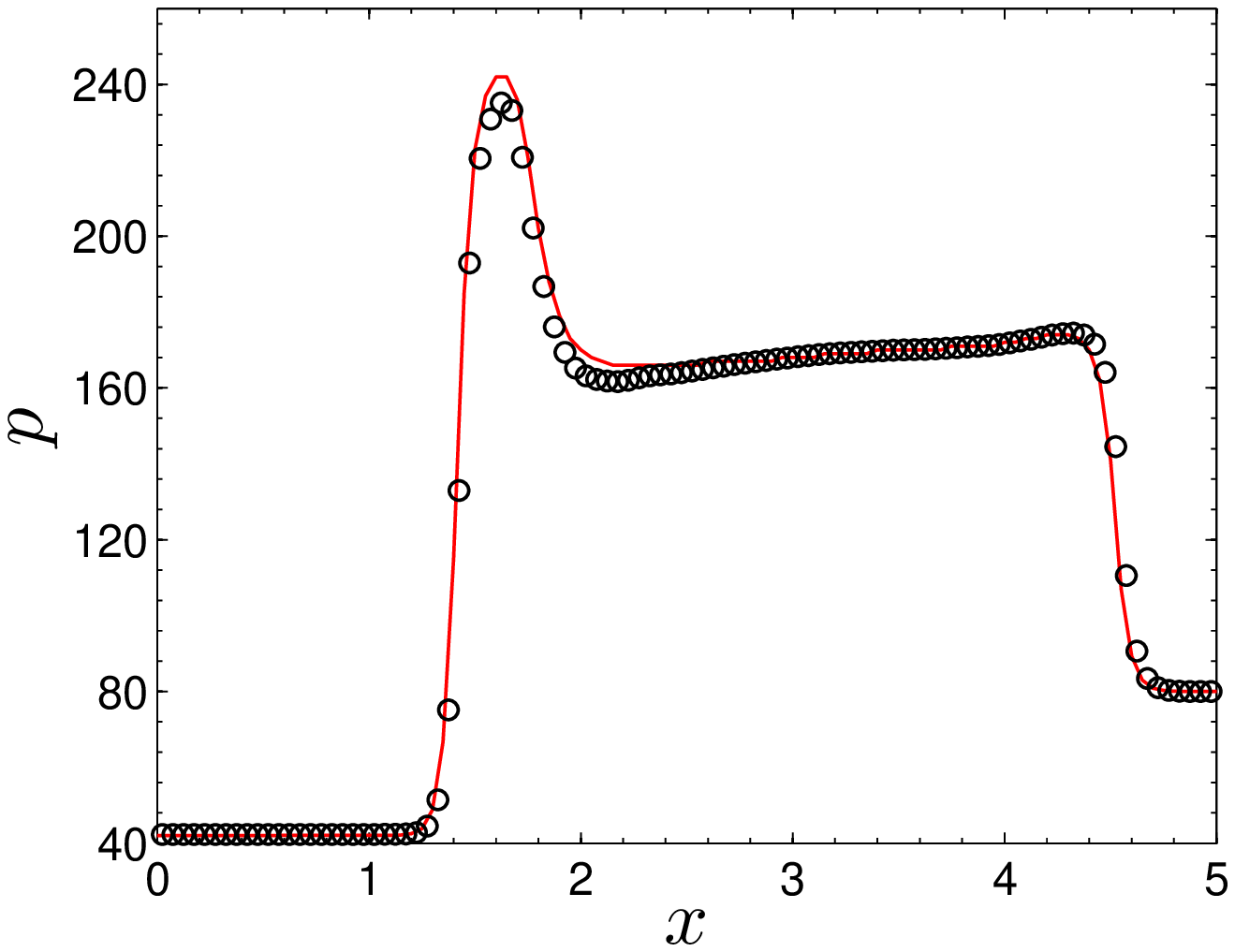}
\includegraphics[width=0.45\textwidth]{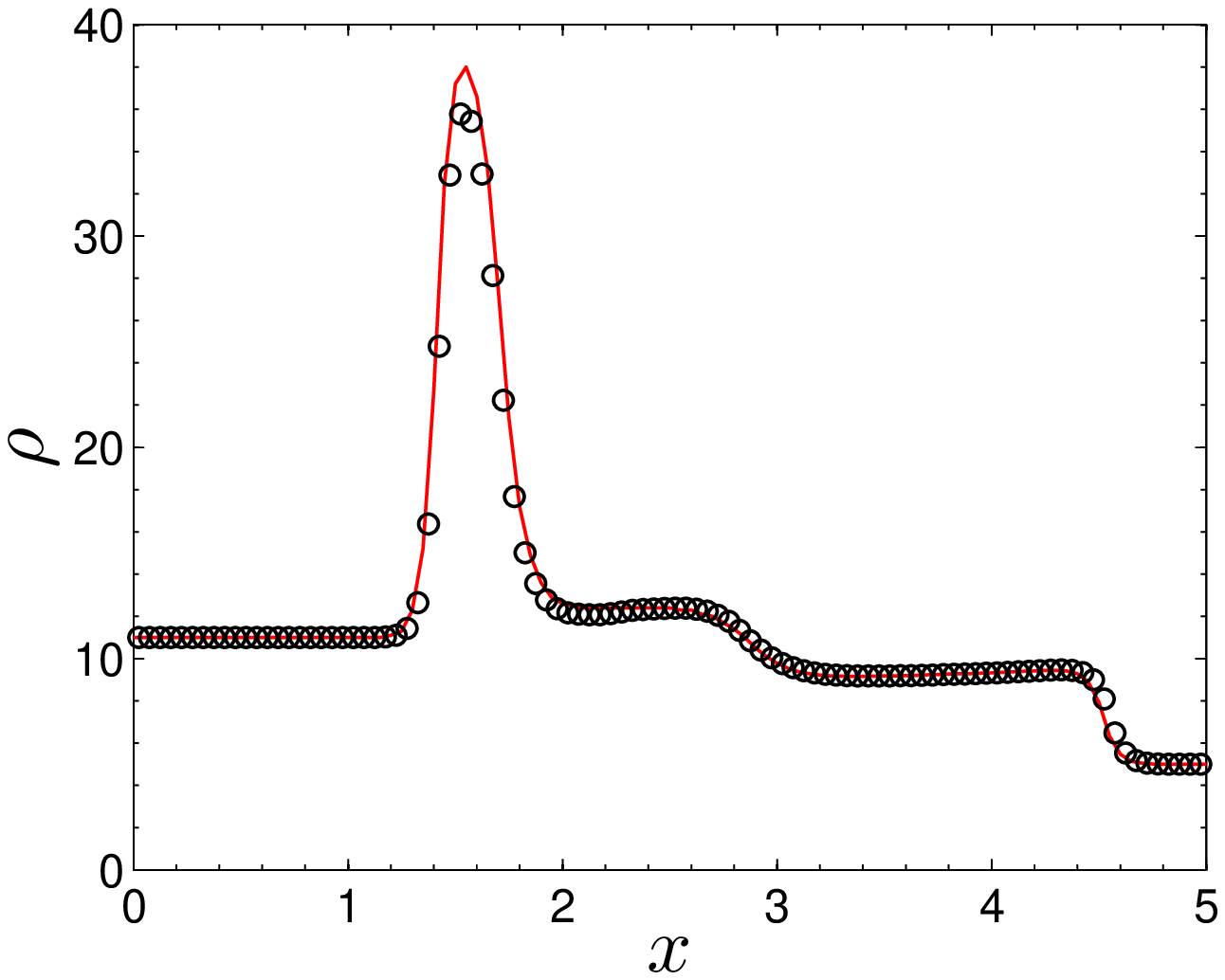}
\caption{Example 6: Comparison of the pressure (left) and density (right) of the present method (circle) and the third-order positivity-preserving scheme in \cite{Huang2019jsc} (line) along the horizontal center line $y=2.5$.}
\label{fig:2D_example2_p_rho_centerline}
\end{figure}

The results of the present boundary treatment and the third-order positivity-preserving scheme in \cite{Huang2019jsc} with reflective boundary conditions are compared in Figs.~\ref{fig:2D_example2_p_rho_contour} (contour plot) and \ref{fig:2D_example2_p_rho_centerline} (solutions along $y=2.5$). It is clear that good agreements between the two methods are achieved again.

\begin{figure}[!ht]
\centering
\includegraphics[width=0.45\textwidth]{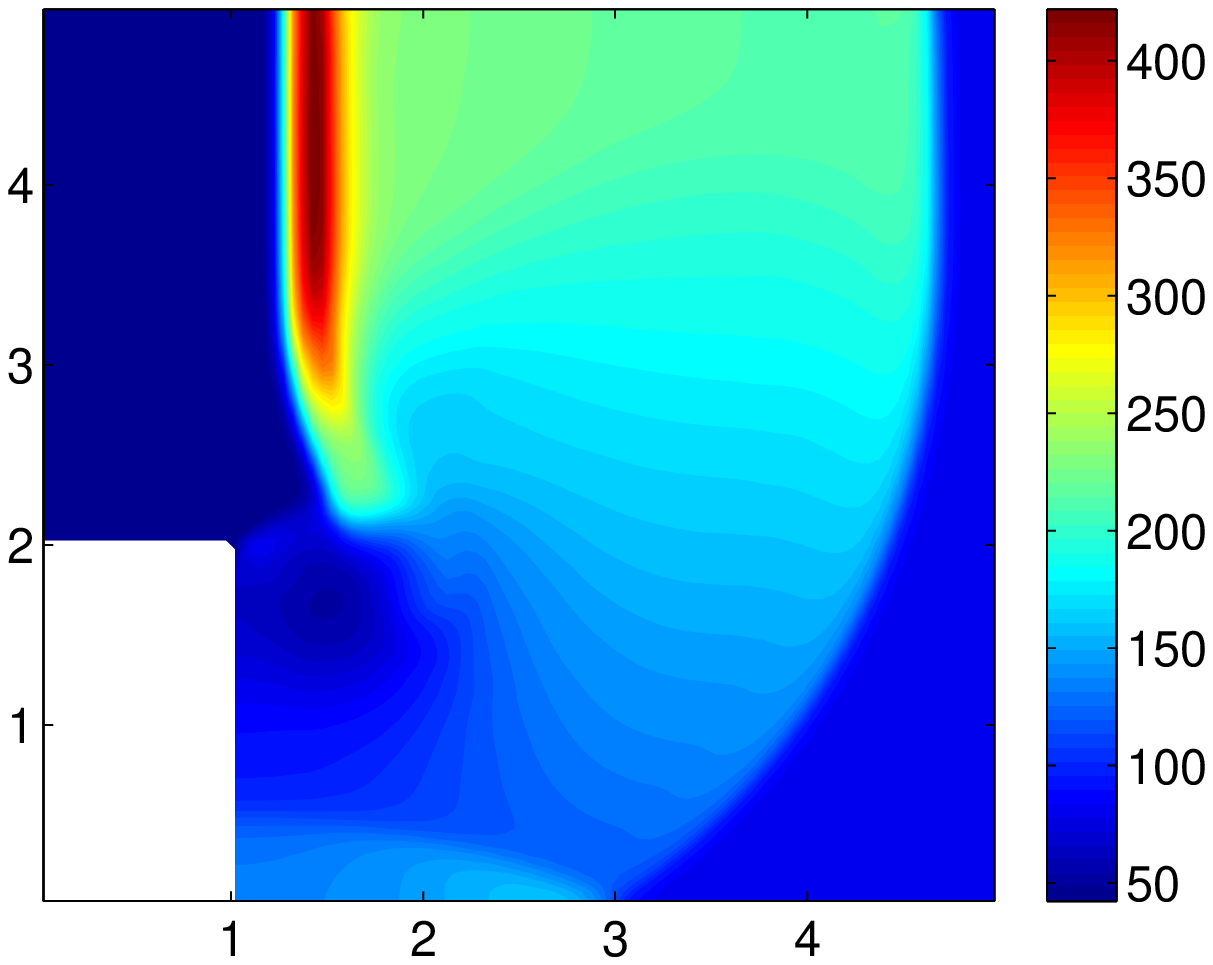}
\includegraphics[width=0.45\textwidth]{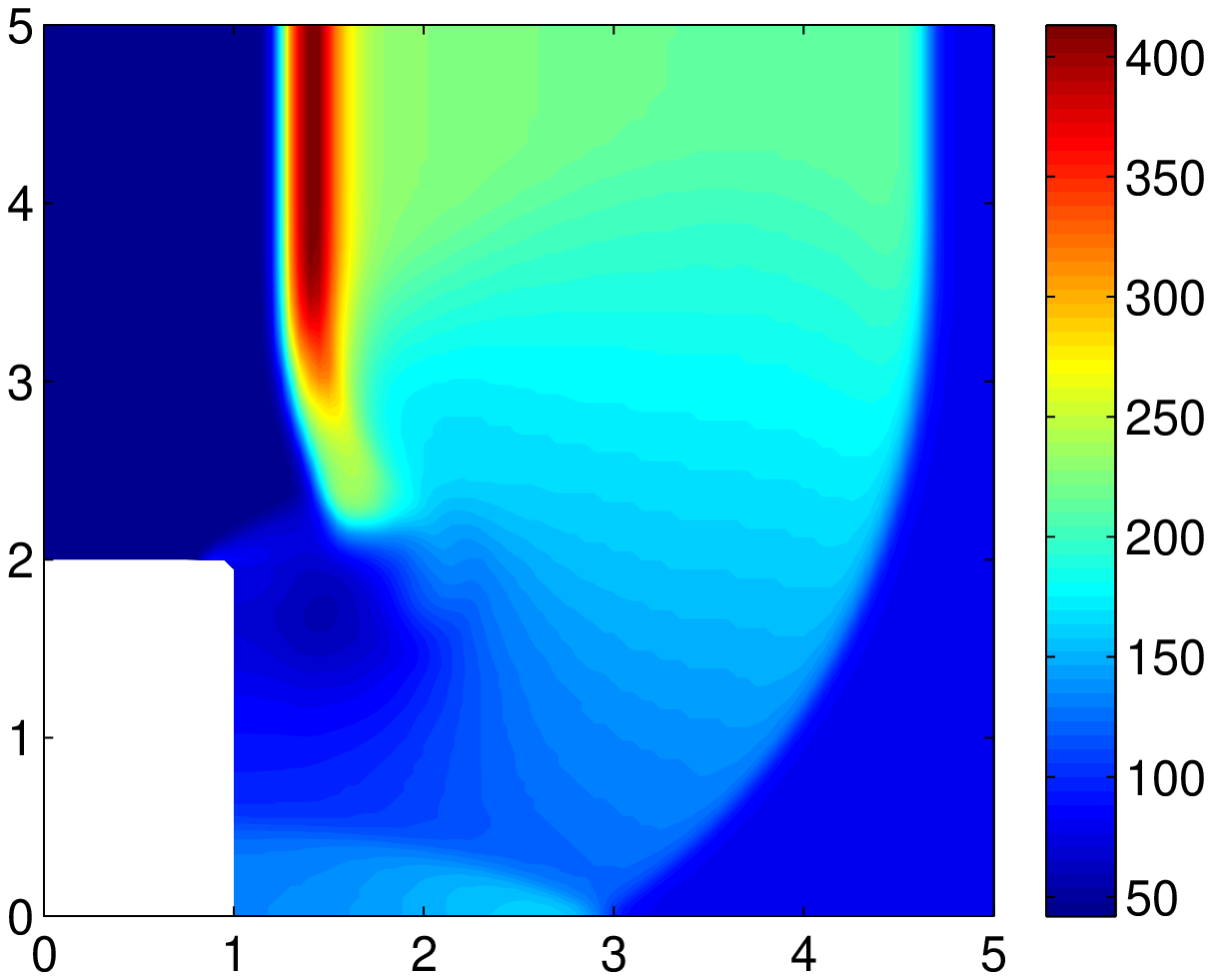}
\includegraphics[width=0.45\textwidth]{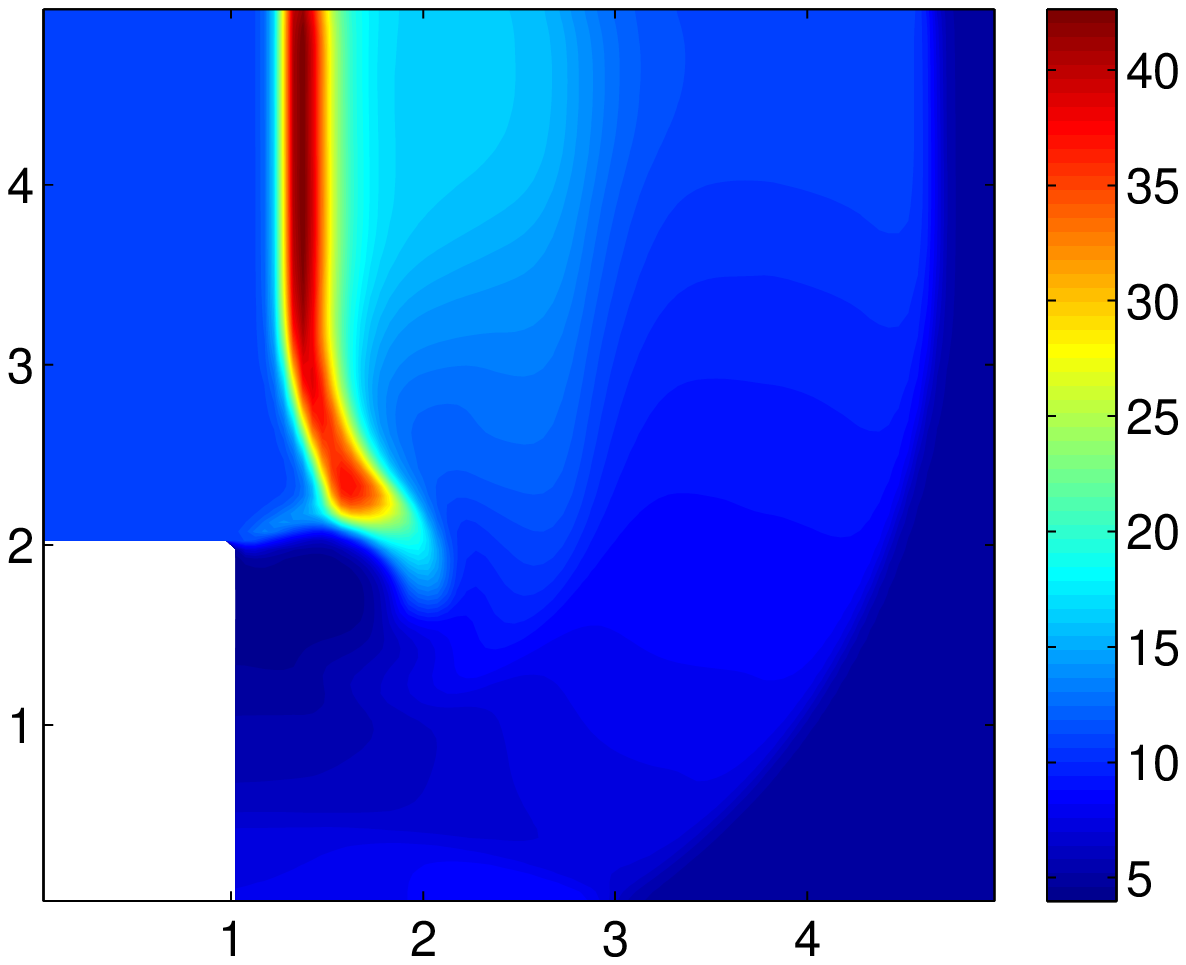}
\includegraphics[width=0.45\textwidth]{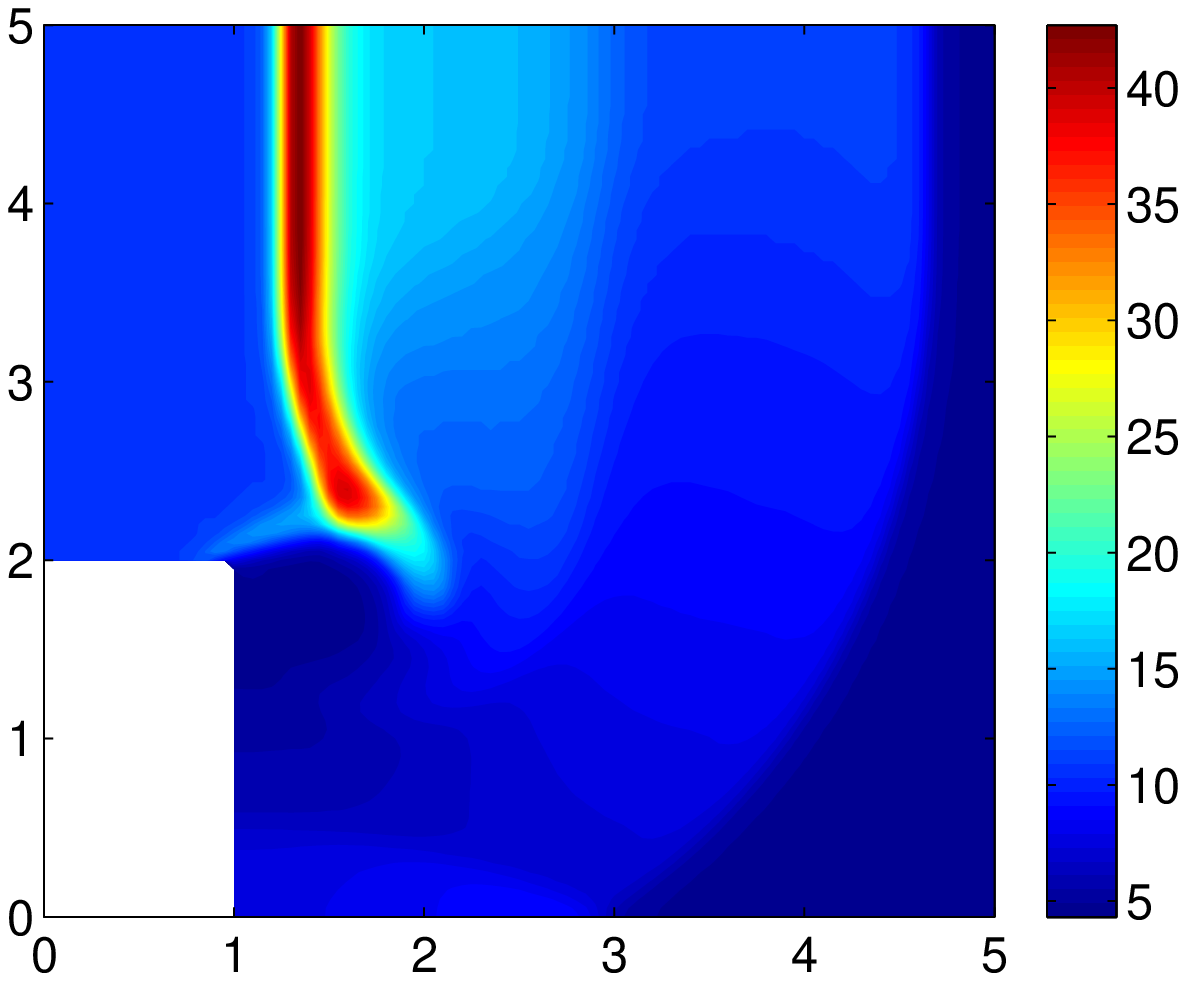}
\caption{Example 6: Contour plots of pressure (top) and density (bottom) of the present method (left) and the third-order positivity-preserving scheme in \cite{Huang2019jsc} (right) with $\Delta x = 1/20$ at $t=0.6$. }
\label{fig:2D_example2_p_rho_contour}
\end{figure}

\begin{figure}[!ht]
\centering
\includegraphics[width=0.45\textwidth]{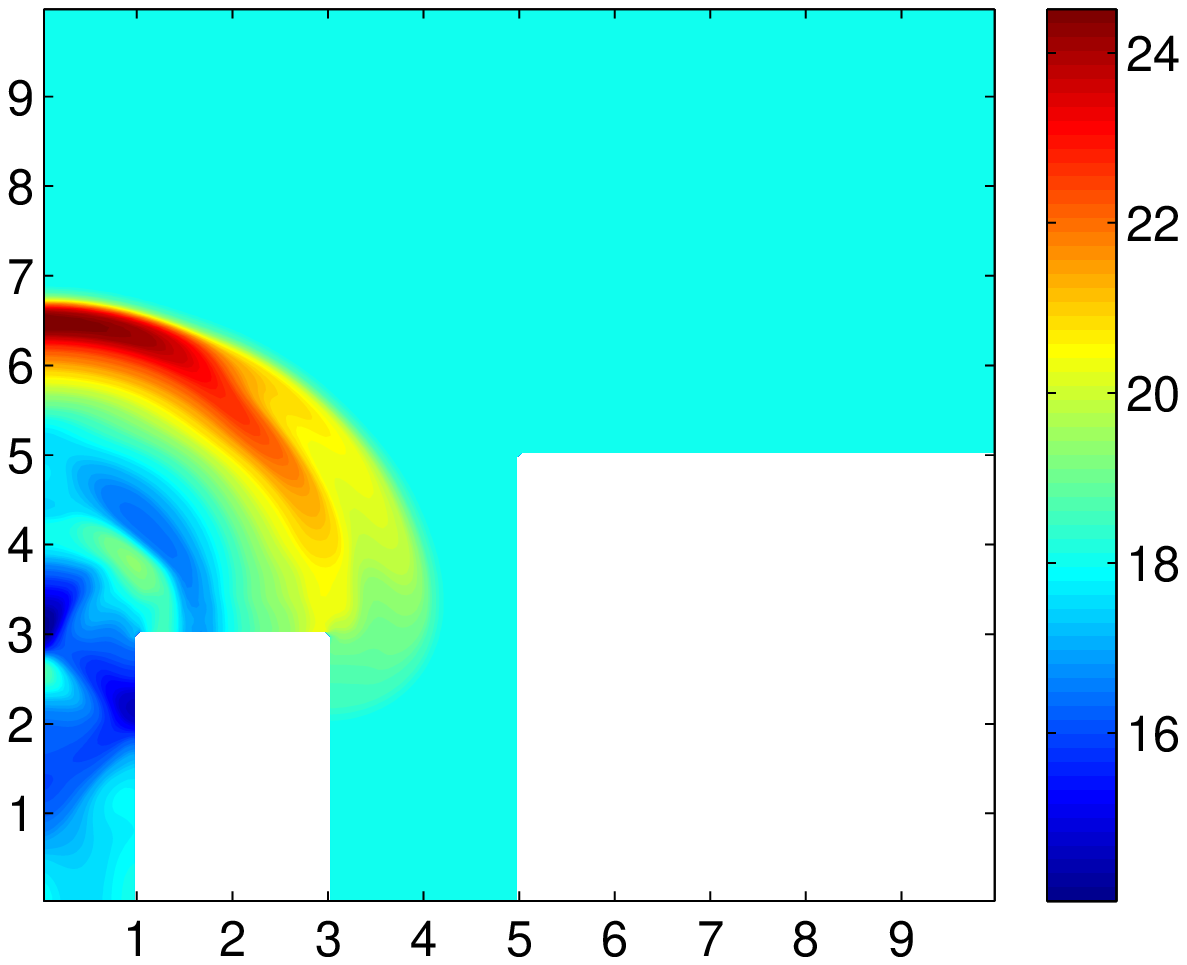}
\includegraphics[width=0.45\textwidth]{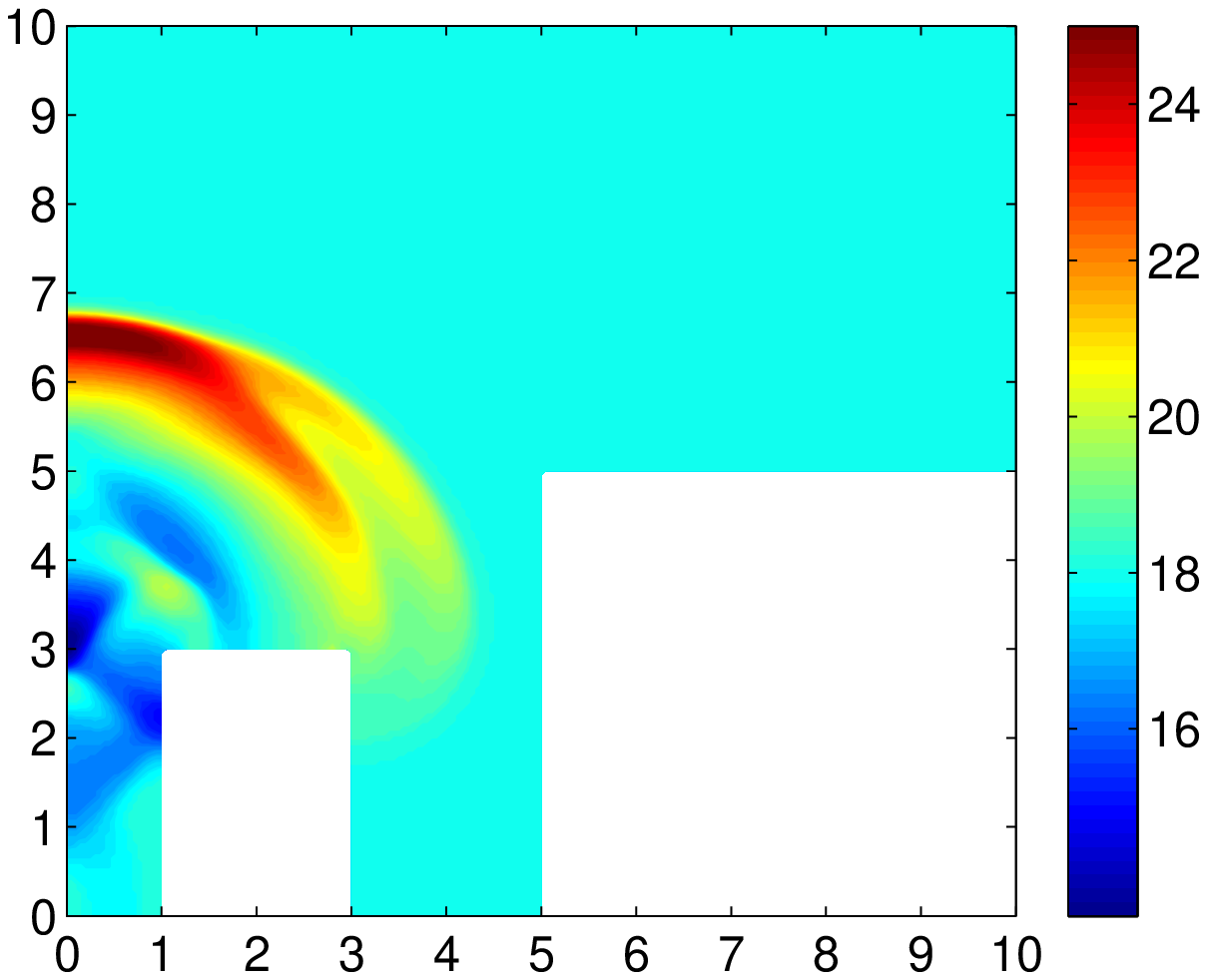}
\includegraphics[width=0.45\textwidth]{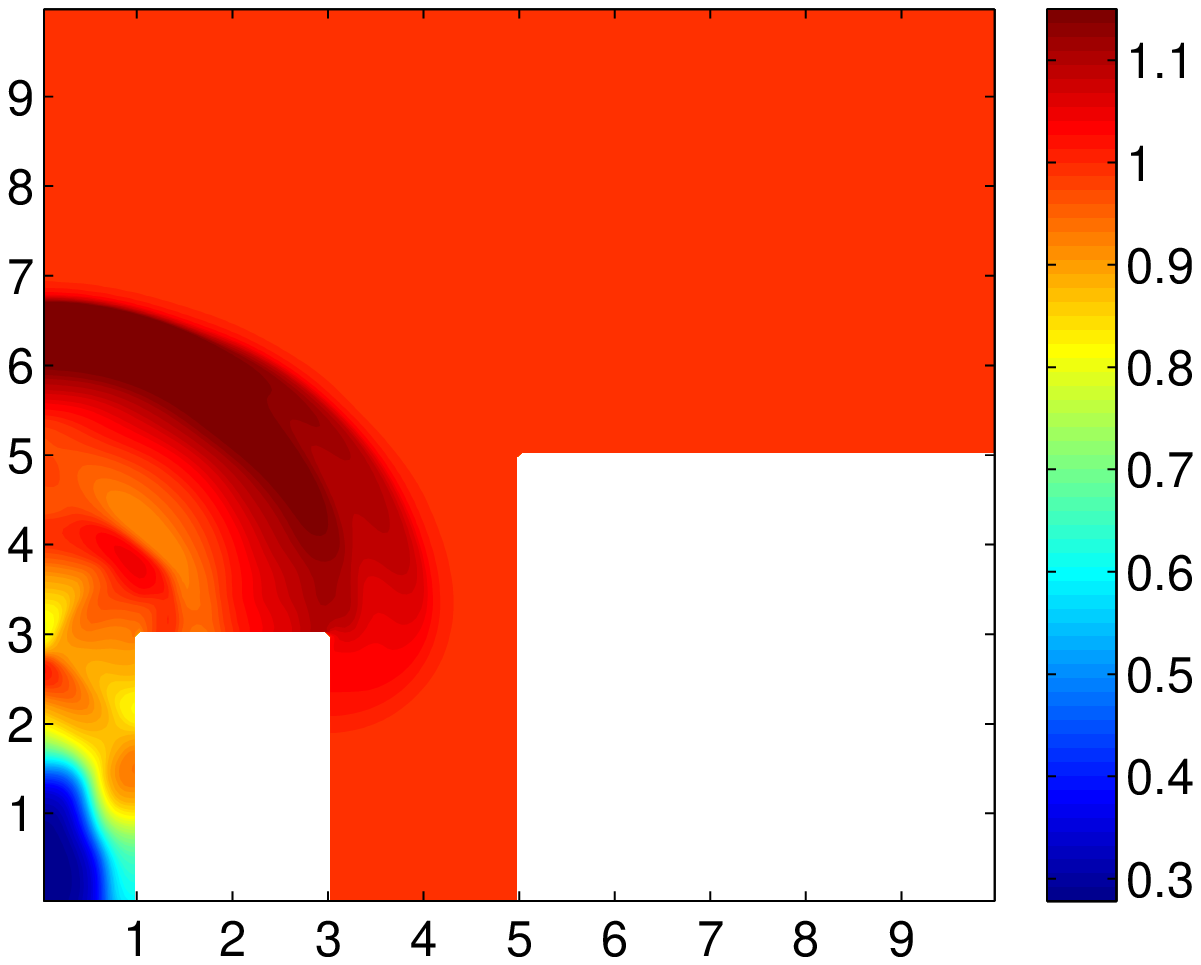}
\includegraphics[width=0.45\textwidth]{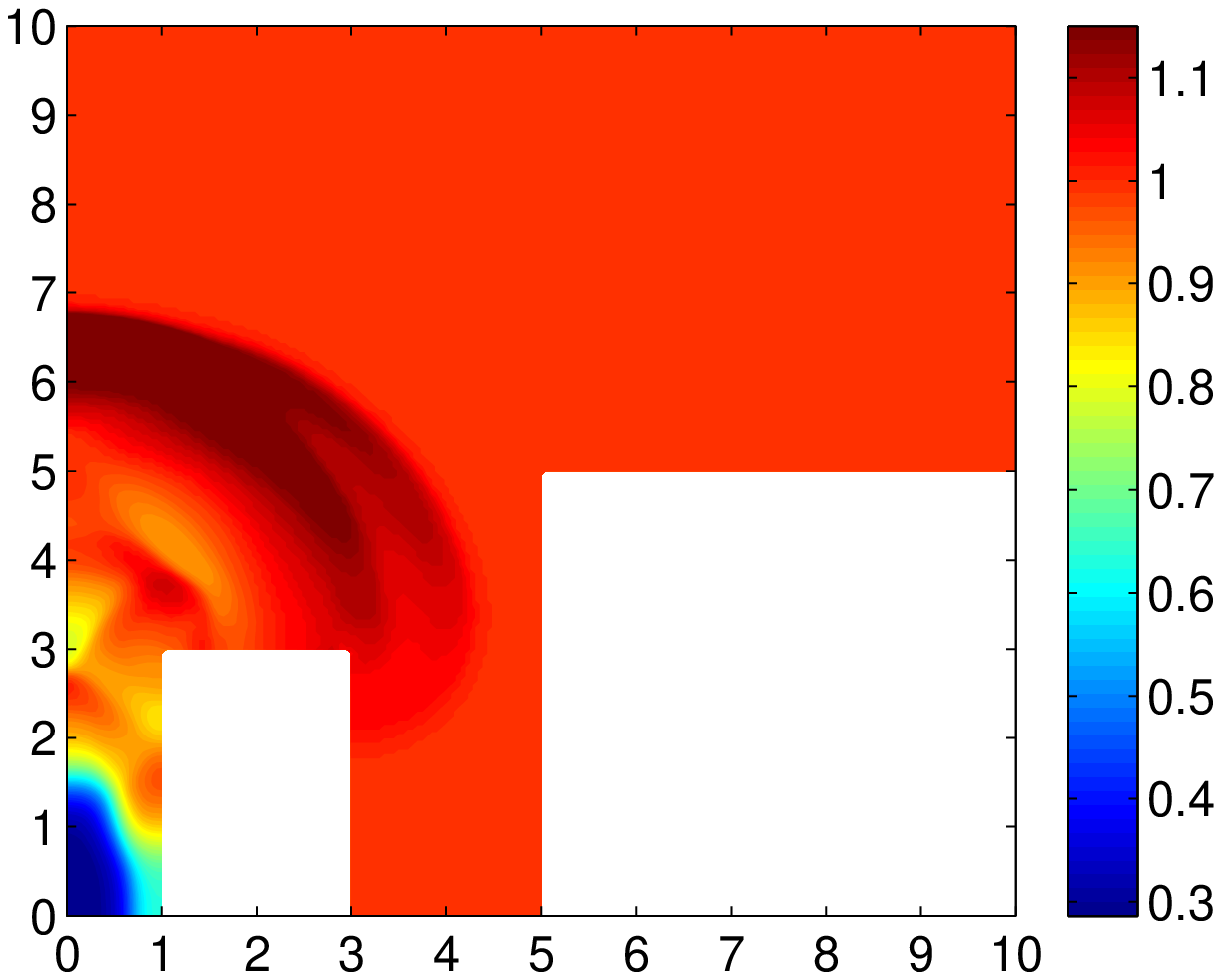}
\caption{Example 7: Contour plots of pressure (top) and density (bottom) of the present method (left) and the third-order positivity-preserving scheme in \cite{Huang2019jsc} (right) with $\Delta x = 1/20$ at $t=1$.}
\label{fig:2D_example3_p_contour}
\end{figure}

\begin{figure}[!ht]
\centering
\includegraphics[width=0.45\textwidth]{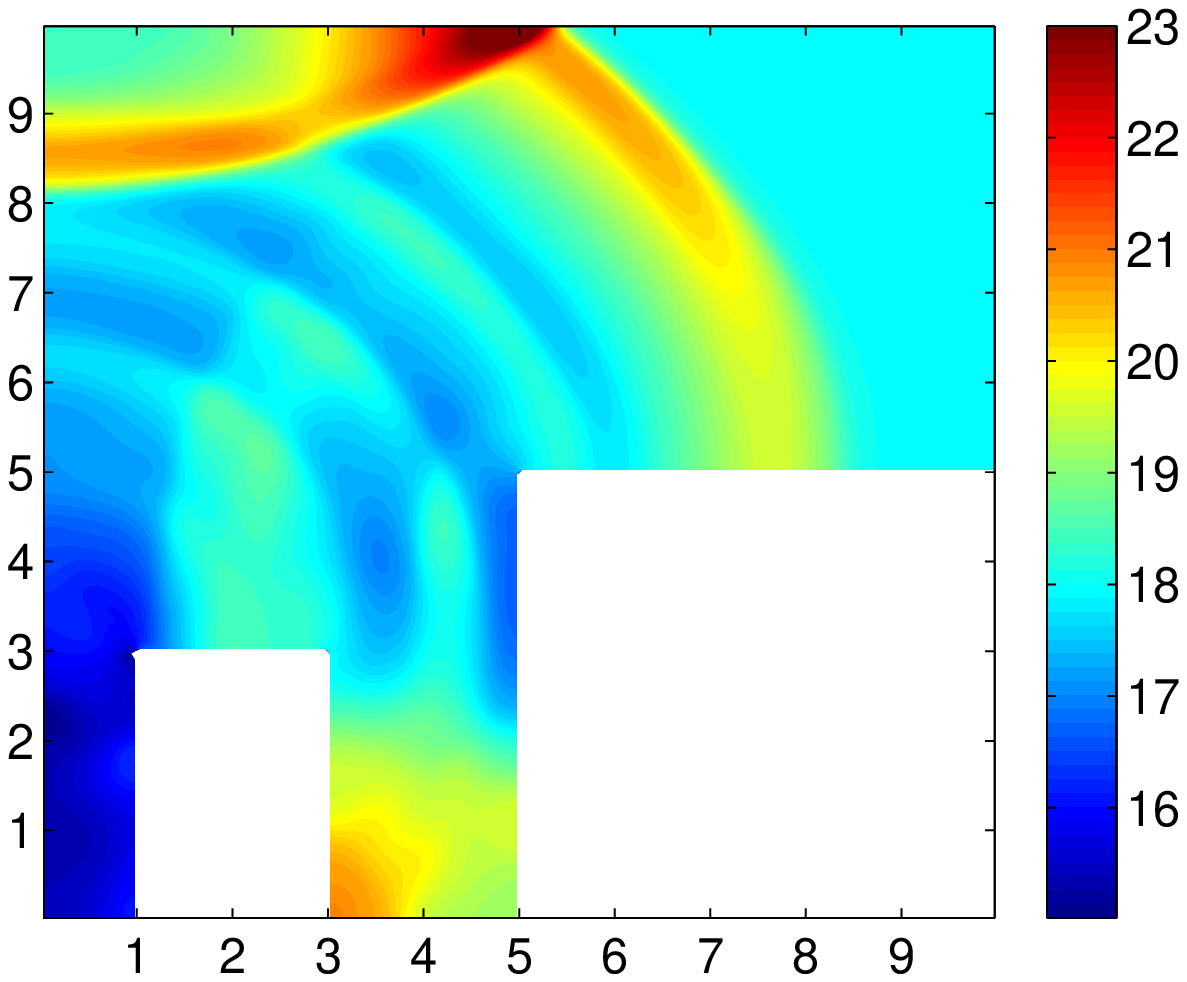}
\includegraphics[width=0.45\textwidth]{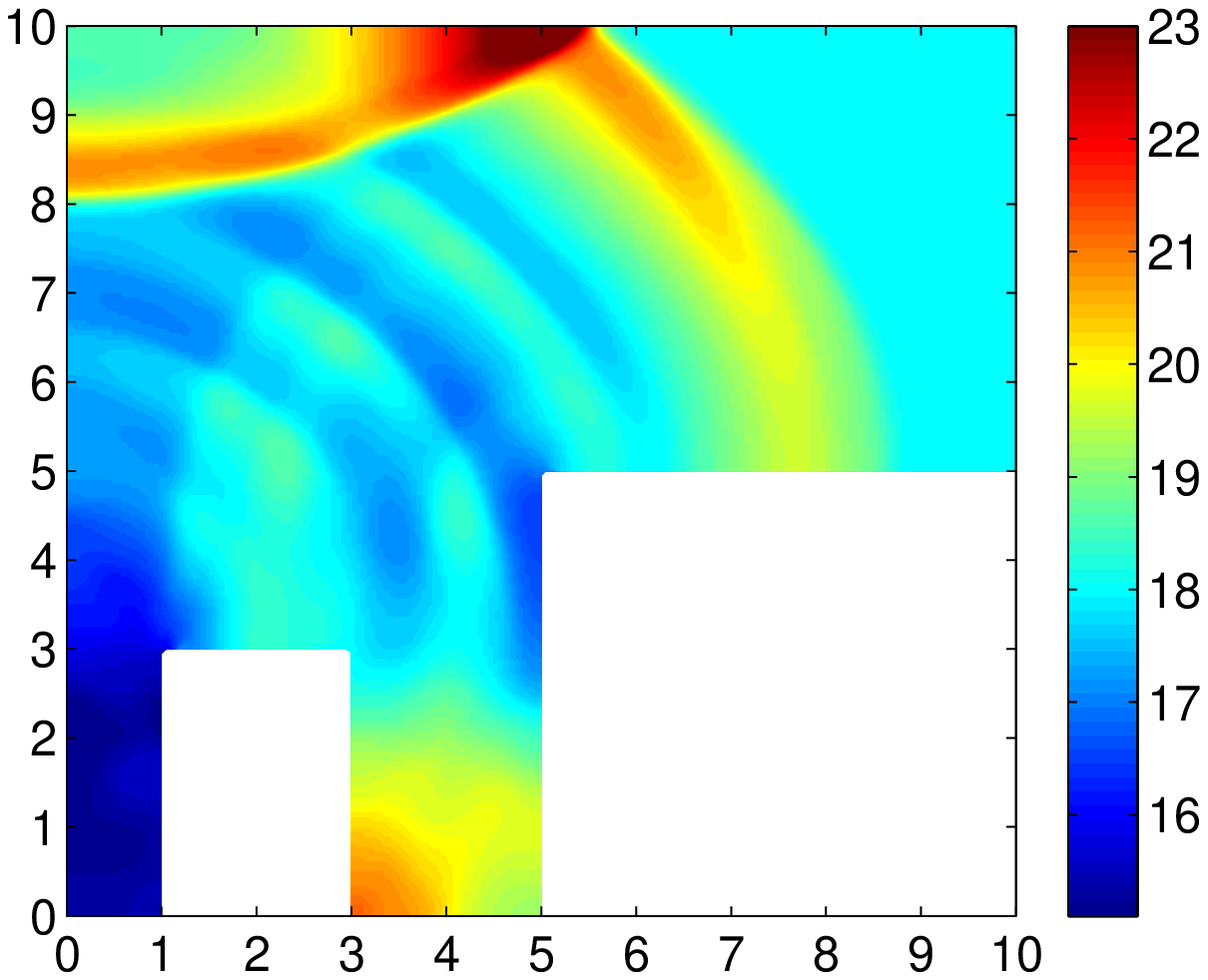}
\includegraphics[width=0.45\textwidth]{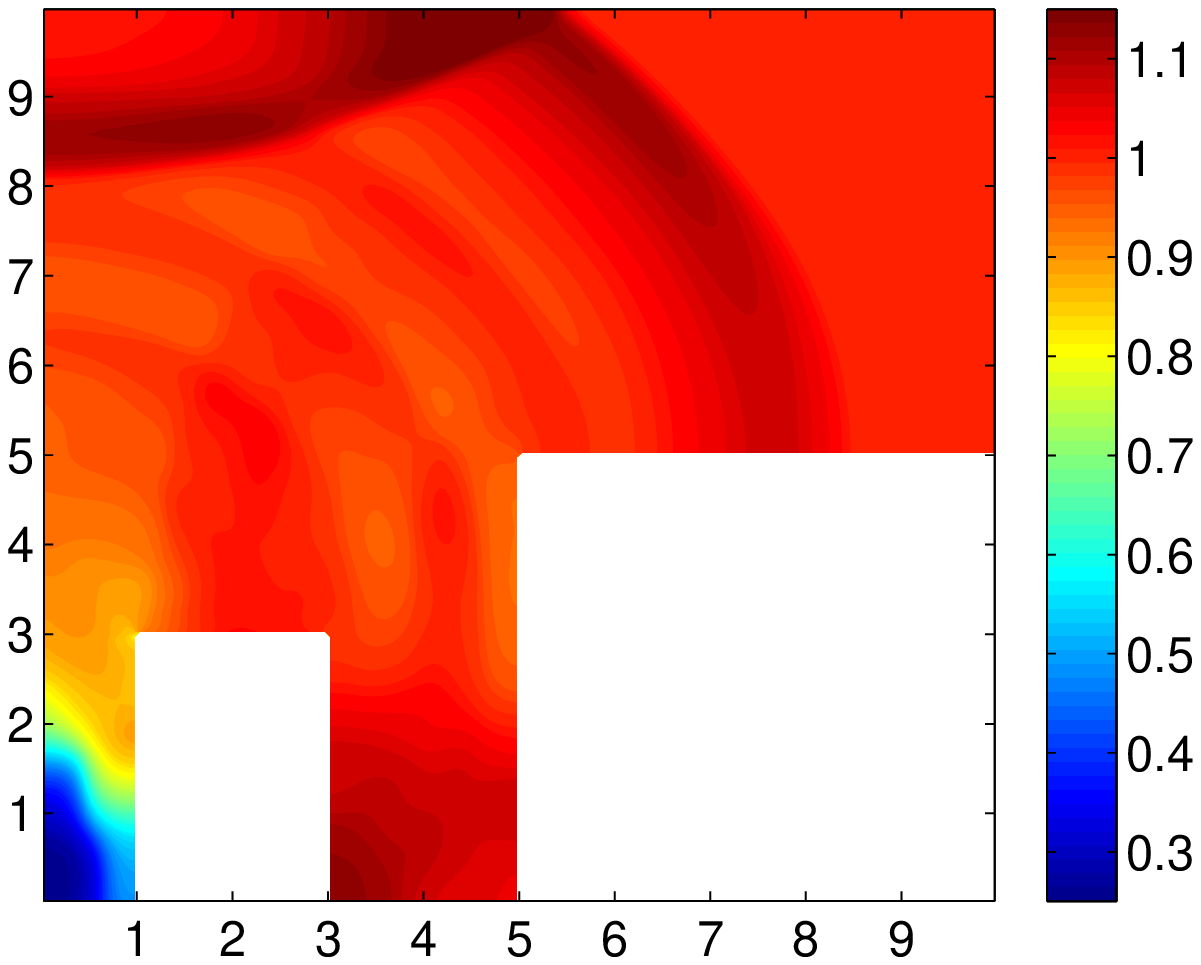}
\includegraphics[width=0.45\textwidth]{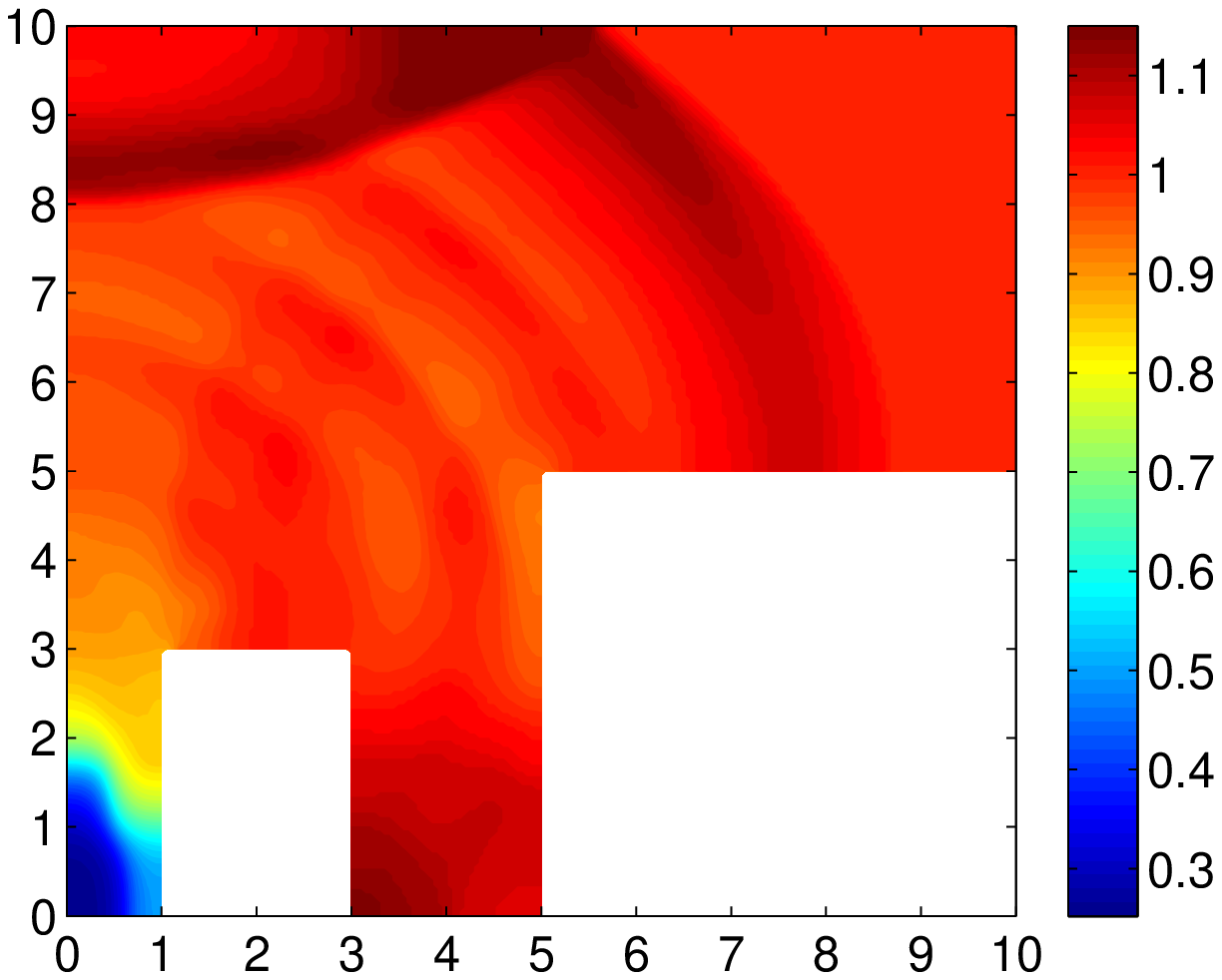}
\caption{Example 7: Contour plots of pressure (top) and density (bottom) of the present method (left) and the third-order positivity-preserving scheme in \cite{Huang2019jsc} (right) with $\Delta x = 1/20$ at $t=2$. }
\label{fig:2D_example3_rho_contour}
\end{figure}

\noindent \textbf{Example 7.}
Finally, we test our method through an example with multiple obstacles in \cite{Huang2019jsc}, which is designed following \cite{Wang2012JCP}.
In this example, the spatial domain is $[0,10] \times [0,10]$ and there are two obstacles with positions $[1,3] \times [0,3]$ and $[5,10]\times [0,5]$, respectively. The initial condition is set as \eqref{2D_example0_InitialCondition} and the boundary condition is $\bm n \cdot \bm u = 0$ for all the boundaries. We take $\Delta x = 1/20$ and results of the present method and the third-order positivity-preserving scheme in \cite{Huang2019jsc} with reflective boundary conditions are shown in Fig.~\ref{fig:2D_example3_p_contour} at $t=1$ and Fig.~\ref{fig:2D_example3_rho_contour} at $t=2$.
It can be seen that the results of the two methods are very similar even with the very coarse mesh, which demonstrate the effectiveness of our boundary treatment for complex domains.

\section{Conclusions and remarks}

In this paper, we propose a high order finite difference boundary treatment method for the IMEX RK schemes solving hyperbolic systems with possibly stiff source terms on a Cartesian mesh. By combining the idea of using the RK schemes at the boundary and the ILW procedure, our method not only preserves the accuracy of the RK schemes but also processes good stability. Our method is different from the widely used approach for the explicit RK schemes by imposing boundary conditions for intermediate solutions \cite{Carpenter1995}, which could not be derived for the IMEX schemes. In addition, the intermediate boundary conditions are only available for explicit RK schemes up to third order while our method applies to arbitrary order IMEX and explicit RK schemes. Moreover, we show that our method applied to the third-order IMEX RK scheme solving 1D scalar equations  can yield the widely used intermediate-stage boundary conditions in \cite{Carpenter1995}. Finally, both 1D examples and 2D reactive Euler equations are used to demonstrate the good stability and third-order accuracy of our boundary treatment for a specific third-order IMEX scheme.

Since the IMEX RK schemes include explicit schemes as special cases, the present boundary treatment also applies to explicit RK schemes. It does not rely on intermediate-stage boundary conditions and thus is expected to be valid for explicit RK schemes of order higher than three. This will be investigated in our next work. Furthermore, our method may be adapted to IMEX RK schemes solving the other partial differential equations, \eg  convection-diffusion equations \cite{Wang2015,lu2016diffusion,Wang2018} and the Boltzmann equation \cite{FY2013,Hu2017,Hu2018,Jin2018}. Many unresolved issues for the boundary treatment of RK schemes can be explored based on the idea of this work.

%
%
%
%


\end{document}